\documentclass[11pt,reqno]{amsart}

\title[Variational and Geometric Structures of Discrete Dirac Mechanics]{Variational and Geometric Structures of\\ Discrete Dirac Mechanics}
\author{Melvin Leok}
\author{Tomoki Ohsawa}
\address{Department of Mathematics, University of California, San Diego, 9500 Gilman Drive, La Jolla, California, USA.}
\email{mleok@math.ucsd.edu, tohsawa@ucsd.edu}
\dedicatory{Communicated by Arieh Iserles}

\date{\today}

\usepackage[margin=1in]{geometry}
\usepackage{amsmath,amsthm,amssymb,graphicx,paralist}
\usepackage[authoryear,sort&compress]{natbib}
\bibpunct{[}{]}{,}{n}{,}{,}

\usepackage[colorlinks=false,bookmarksdepth=0]{hyperref}
\usepackage[all]{xy}
\entrymodifiers={+!!<0pt,\fontdimen22\textfont2>}

\theoremstyle{plain}
\newtheorem{theorem}{Theorem}[section]
\newtheorem{proposition}[theorem]{Proposition}

\theoremstyle{definition}
\newtheorem{definition}[theorem]{Definition}
\newtheorem{example}[theorem]{Example}

\theoremstyle{remark}
\newtheorem{remark}[theorem]{Remark}

\numberwithin{equation}{section}

\def\FL{\mathbb{F}L}
\def\defeq{\mathrel{\mathop:}=}
\def\setdef#1#2{ \left\{ #1 \ |\ #2 \right\} }
\def\R{\mathbb{R}}
\def\pd#1#2{\dfrac{\partial #1}{\partial #2}}
\def\pdd#1#2#3{\dfrac{\partial^2#1}{\partial #2\partial #3}}
\def\tpd#1#2{\partial #1/\partial #2}

\def\parentheses#1{\left(#1\right)}
\def\braces#1{\left\{#1\right\}}
\def\brackets#1{\left[#1\right]}
\def\rank{\mathop{\mathrm{rank}}\nolimits}

\def\Span{\mathop{\mathrm{span}}\nolimits} % span << WARNING: Defining "span" makes gather and align impossible to be used. Hence "Span."

\def\ip#1#2{\left\langle#1,#2\right\rangle}
\def\DS{\displaystyle}
\def\id{\mathop{\mathrm{id}}\nolimits}

\begin{document}
\footskip=.6in

%\allowdisplaybreaks

\begin{abstract}
In this paper, we develop the theoretical foundations of discrete Dirac mechanics, that is, discrete mechanics of degenerate Lagrangian/Hamiltonian systems with constraints. We first construct discrete analogues of Tulczyjew's triple and induced Dirac structures by considering the geometry of symplectic maps and their associated generating functions. We demonstrate that this framework provides a means of deriving discrete Lagrange--Dirac and nonholonomic Hamiltonian systems.
In particular, this yields nonholonomic Lagrangian and Hamiltonian integrators. We also introduce discrete Lagrange--d'Alembert--Pontryagin and Hamilton--d'Alembert variational principles, which provide an alternative derivation of the same set of integration algorithms. The paper provides a unified treatment of discrete Lagrangian and Hamiltonian mechanics in the more general setting of discrete Dirac mechanics, as well as a generalization of symplectic and Poisson integrators to the broader category of Dirac integrators.
\end{abstract}

\maketitle

\begin{center}
{\it Dedicated to the memory of Jerrold E. Marsden.}
\end{center}

\section{Introduction}
Dirac structures, which can be viewed as simultaneous generalizations of symplectic and Poisson structures, were introduced in \citet{Co1990a,Co1990b}. In the context of geometric mechanics~\citep{Ar1989, AbMa1978, MaRa1999}, Dirac structures are of interest as they can directly incorporate Dirac constraints that arise in degenerate Lagrangian systems~\cite{Di1950,Di1958a,Di1964,Ku1969,GoNe1979a,GoNe1979b,GoNe1980}, interconnected systems~\citep{CeScBa2003, Sc2006}, and nonholonomic systems~\citep{Bl2003}, and thereby provide a unified geometric framework for studying such problems.

From the Hamiltonian perspective, these systems are described by implicit Hamiltonian systems; see \citet{BlCr1997} and \citet{Sc1998} for applications of such a formulation to LC circuits and nonholonomic systems, and \citet{DaSc1998} for a comprehensive review of Dirac structures in this setting.
This approach is motivated by earlier work on almost-Poisson structures that describe nonholonomic systems using brackets that fail to satisfy the Jacobi identity~\citep{ScMa1994}.
These ideas are further extended to define port-Hamiltonian systems, which are intended to model interconnected systems~(see \citet{Sc2006} for a survey of such applications).
 
On the Lagrangian side,  degenerate, interconnected, and nonholonomic systems can be described by Lagrange--Dirac (or implicit Lagrangian) systems introduced by \citet{YoMa2006a} in the context of Tulczyjew's triple~\citep{Tu1976a, Tu1976b} and a certain class of representations of Dirac structures called induced Dirac structures~\citep{DaSc1998}.
The resulting Lagrange--Dirac equations generalize the Lagrange--d'Alembert equations for nonholonomic systems.
The corresponding variational description of Lagrange--Dirac systems was developed in \citet{YoMa2006b}, with the introduction of the Hamilton--Pontryagin principle on the Pontryagin bundle $TQ \oplus T^*Q$, which yields the generalized Legendre transformation, as well as Hamilton's principle for Lagrangian systems and Hamilton's phase space principle for Hamiltonian systems.
\citet{YoMa2006b} also introduced the Lagrange--d'Alembert--Pontryagin principle, a generalization of the Hamilton--Pontryagin principle, which yields Lagrange--Dirac systems with nonholonomic constraints.
It also generalizes the Lagrange--d'Alembert principle for nonholonomic systems~(see, e.g., \citet{Bl2003}).

In the context of geometric numerical integration~\citep{HaLuWa2006, LeRe2004}, which is concerned with the development of numerical methods that preserve geometric properties of the corresponding continuous flow, variational integrators that preserve the symplectic structure can be systematically derived from a discrete Hamilton's principle~\citep{MaWe2001}, and can be extended to asynchronous variational integrators~\citep{LeMaOrWe2003} that preserve the multisymplectic structure of Hamiltonian partial differential equations. The discrete variational formulation of Hamiltonian mechanics was developed by \citet{LaWe2006} as the dual, in the sense of optimization, to discrete Lagrangian mechanics. Discrete analogues of the Hamilton--Pontryagin principle were introduced in \cite{Kh2006, BoMa2008} for particular choices of discrete Lagrangians.
Discrete Lagrangian, Hamiltonian, and nonholonomic mechanics have also been generalized to Lie groupoids~\citep{We1996a, MaMaMa2006, IgMaMaMa2008, St2010}.

\subsection*{Contributions of this paper}
In this paper, we introduce discrete analogues of Tulczyjew's triple and induced Dirac structures, and show how they describe discrete Lagrange--Dirac and nonholonomic Hamiltonian systems.
The construction relies on the observation that Tulczyjew's triple arises from symplectic maps between the iterated tangent and cotangent bundles $T^{*}TQ$, $TT^{*}Q$, and $T^{*}T^{*}Q$.
By analogy, we construct discrete analogues of Tulczyjew's triple that are derived from properties of symplectic maps between discrete analogues of the iterated tangent and cotangent bundles.
We then demonstrate that they yield discrete Lagrange--Dirac and nonholonomic Hamiltonian systems, and recover nonholonomic integrators that are typically derived from a discrete Lagrange--d'Alembert principle.

We also introduce discrete analogues of the Lagrange--d'Alembert--Pontryagin and Hamilton--d'Alembert variational principles, which provide a variational characterization of discrete Lagrange--Dirac and nonholonomic Hamiltonian systems that we previously described in terms of the discrete analogues of Tulczyjew's triple and induced Dirac structures.
The discrete Lagrange--Dirac and nonholonomic Hamiltonian systems recover the standard Lagrangian variational integrators (see, e.g., \citet{MaWe2001}), Hamiltonian variational  integrators of \citet{LaWe2006}, and nonholonomic integrators~(see, e.g., \citet{CoMa2001} and \citet{McPe2006}).

Discrete Hamiltonian mechanics \citep{LaWe2006} is not intrinsic, due to its dependence on Type~2 or 3 generating functions of symplectic maps. Since discrete Dirac mechanics encompasses discrete Hamiltonian mechanics, we first limit our discussions to the cases where the configuration manifold $Q$ is a vector space.
We then introduce a retraction, a map from $TQ$ to $Q$, to extend the ideas to the more general case where $Q$ is a manifold.
Specifically, we extend the Lagrange--d'Alembert--Pontryagin principle to this case, and show that it yields, using a certain class of coordinate charts specified by the retraction, the same coordinate expressions for Lagrange--Dirac systems as in the linear case.
This gives a firm theoretical foundation and a prescription for performing computations with Lagrange--Dirac systems on manifolds.

\subsection*{Outline of this paper}
The paper is organized as follows.
In Section~\ref{sec:DS_TT_LDS}, we review induced Dirac structures, Tulczyjew's triple, and Lagrange--Dirac systems with an LC circuit as a motivating example.
In Sections~\ref{sec:DATT} and \ref{sec:DAIDS}, we construct discrete analogues of Tulczyjew's triple and induced Dirac structures.
These discrete analogues lead us to the development of discrete Dirac mechanics, i.e., discrete Lagrange--Dirac and nonholonomic Hamiltonian systems, in Section~\ref{sec:DDM}.
We then come back to the LC circuit example in Section~\ref{sec:LC_circuit}: We discretize the LC circuit and describe it as a discrete Lagrange--Dirac system to obtain a numerical method; we also test the method numerically and compare the result with an exact solution.
In Section~\ref{sec:VariationalStructure}, we briefly come back to the continuous-time setting to review the Lagrange--d'Alembert--Pontryagin and Hamilton--d'Alembert principles for Lagrange--Dirac and nonholonomic Hamiltonian systems.
Then, in Section~\ref{sec:DiscreteVariationalStructure}, we define the discrete analogues of the variational principles.
In Section~\ref{sec:ExtensionToManifolds}, we extend our results to computations on manifolds.

\section{Dirac Structures, Tulczyjew's Triple, and Lagrange--Dirac Systems}
\label{sec:DS_TT_LDS}
We first briefly review the induced Dirac structures that give rise to Lagrange--Dirac systems, taking an LC circuit as an example~(see \cite{YoMa2006a, YoMa2006b, YoMa2007b}).
Lagrange--Dirac systems are particularly useful in formulating systems with degenerate Lagrangians and/or constraints.
LC circuits are a class of examples that is particularly well suited for the formulation as Lagrange--Dirac systems, since they often involve degenerate Lagrangians and also constraints arising from the Kirchhoff laws.

\subsection{LC Circuit---Example of Degenerate Lagrangian System with Constraints}
\begin{figure}[htbp]
  \centering
  \includegraphics[width=.375\linewidth]{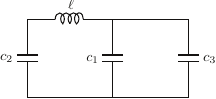}
  \caption{LC circuit---Example of degenerate Lagrangian system with constraints (see \citep{YoMa2006a}).}
  \label{fig:LCC}
\end{figure}

Following \citet{YoMa2006a}, consider the LC circuit with an inductor $\ell$ and three capacitors $c_{1}, c_{2}$, and $c_{3}$ shown in Fig.~\ref{fig:LCC}.
The configuration space is the 4-dimensional vector space $Q = \{(q^{\ell}, q^{c_{1}}, q^{c_{2}}, q^{c_{3}})\}$, which represents charges in the circuit elements.
Then, an element $f_{q} = (f^{\ell}, f^{c_{1}}, f^{c_{2}}, f^{c_{3}})$ in the tangent space $T_{q}Q$ represents the currents in the corresponding circuit elements; hence the tangent bundle $TQ$ is a charge-current space.
The Lagrangian $L: TQ \to \R$ is given by
\begin{equation}
  \label{eq:L-LC}
  L(q,f) = \frac{\ell}{2} (f^{\ell})^{2} - \frac{(q^{c_{1}})^{2}}{2c_{1}} - \frac{(q^{c_{2}})^{2}}{2c_{2}} - \frac{(q^{c_{3}})^{2}}{2c_{3}}.
\end{equation}
The Lagrangian is clearly degenerate:
\begin{equation*}
  \det \parentheses{ \pdd{L}{f^{i}}{f^{j}} } = 0,
\end{equation*}
which corresponds to the fact that not every circuit component has inductance.
Therefore, the Legendre transformation $\FL: TQ \to T^{*}Q$,  with $T^{*}Q$ being the cotangent bundle of $Q$, defined by
\begin{equation*}
  \FL: f \mapsto \pd{L}{f^{i}} dq^{i}
\end{equation*}
is not invertible, and hence it is impossible to write the system as a Hamiltonian system in the conventional sense.
Notice also that the Kirchhoff current law imposes the constraints $-f^{\ell} + f^{c_{2}} = 0$ and $-f^{c_{1}} + f^{c_{2}} - f^{c_{3}} = 0$.
This defines the constraint distribution $\Delta_{Q} \subset TQ$ given by
\begin{equation}
  \label{eq:KCL}
  \Delta_{Q} = \setdef{ f \in TQ }{ \omega^{a}(f) = 0,\, a = 1,2 },
\end{equation}
with the constraint one-forms $\{ \omega^{1}, \omega^{2} \}$ defined as
\begin{equation}
  \label{eq:omega-LC}
  \omega^{1} = -dq^{\ell} + dq^{c_{2}},
  \qquad
  \omega^{2} = -dq^{c_{1}} + dq^{c_{2}} - dq^{c_{3}}.
\end{equation}
Then, one can write the constraints simply as $f \in \Delta_{Q}$.
If we introduce the annihilator distribution (or codistribution) $\Delta_{Q}^{\circ} \subset T^{*}Q$ of $\Delta_{Q} \subset TQ$ by
\begin{equation}
  \label{eq:Delta_Q^circ}
  \Delta_{Q}^\circ(q) \defeq
  \setdef{ \alpha_q\in T^*_q Q }{ \forall v_q\in\Delta_{Q},\; \langle \alpha_q, v_q\rangle = 0 },
\end{equation}
then we have $\Delta_{Q}^{\circ} = \Span\{ \omega^{1}, \omega^{2} \}$.

\subsection{Induced Dirac Structures}
\label{ssec:Induced_Dirac_Structure}
The key idea in formulating Lagrange--Dirac systems for systems with constraints like the above LC circuits is to introduce a Dirac structure induced by the above constraints.
Let us first recall the basic definitions and results following \citet{YoMa2006a}.
\begin{definition}[Dirac Structures on Vector Spaces]
  Let $V$ be a vector space and $V^{*}$ be its dual.
  For a subspace $D \subset V \oplus V^{*}$, we define
  \begin{equation}
    D^{\perp} \defeq \setdef{ (v, \alpha) \in V \oplus V^{*} }{ \ip{\alpha'}{v} + \ip{\alpha}{v'} = 0 \text{ for any $(v',\alpha') \in D$} },
  \end{equation}
  where $\ip{\cdot}{\cdot}: V^{*} \times V \to \R$ is the natural pairing.
  A subspace $D$ of $V \oplus V^{*}$ is called a {\em Dirac structure} on $V$ if $D^{\perp} = D$.
\end{definition}
This definition naturally extends to manifolds:
\begin{definition}[Dirac Structures on Manifolds]
  Let $M$ be a manifold and $TM$ and $T^{*}M$ be its tangent and cotangent bundles.
  For a subbundle $D \subset TM \oplus T^{*}M$, we define
  \begin{equation}
    D^{\perp} \defeq \setdef{ (v, \alpha) \in TM \oplus T^{*}M }{ \ip{\alpha'}{v} + \ip{\alpha}{v'} = 0 \text{ for any $(v',\alpha') \in D$} },
  \end{equation}
  where $\oplus$ is the Whitney sum, and $\ip{\cdot}{\cdot}: T^{*}M \times TM \to \R$ is the natural pairing.
  A subbundle $D$ over $M$ of $TM \oplus T^{*}M$ is called a {\em (generalized) Dirac structure} on $M$ if $D^{\perp} = D$.
\end{definition}

A particularly important class of Dirac structures is the {\em induced Dirac structure} on a cotangent bundle defined in the following way:
Let $Q$ be a manifold, $\pi_{Q}: T^{*}Q \to Q$ be the cotangent bundle projection, and $\Omega^{\flat}: TT^{*}Q \to T^{*}T^{*}Q$ be the flat map associated with the standard symplectic structure $\Omega$ on $T^{*}Q$.
\begin{proposition}[The Induced Dirac Structure on $T^{*}Q$; see \citep{YoMa2006a,DaSc1998,Sc1998}]
  \label{prop:IDS}
  Given a constant-dimensional distribution $\Delta_{Q} \subset TQ$ on $Q$, define the lifted distribution
  \begin{equation}
    \label{eq:Delta_T^*Q}
    \Delta_{T^{*}Q} \defeq (T\pi_{Q})^{-1}(\Delta_{Q}) \subset TT^{*}Q,
  \end{equation}
  and let $\Delta_{T^{*}Q}^{\circ} \subset T^{*}T^{*}Q$ be its annihilator, which is also given by $\Delta_{T^{*}Q}^{\circ}=\pi_Q^*(\Delta_Q^\circ)$.
  Then, the subbundle $D_{\Delta_{Q}} \subset TT^{*}Q \oplus T^{*}T^{*}Q$ defined by
  \begin{equation}
    \label{eq:D_Delta_Q}
    D_{\Delta_{Q}} \defeq
    \setdef{ (v,\alpha)\in TT^{*}Q \oplus T^{*}T^{*}Q }
    { v \in \Delta_{T^*Q},
      \;
      \alpha - \Omega^{\flat}(v) \in \Delta_{T^{*}Q}^{\circ} }
  \end{equation}
  is a Dirac structure on $T^{*}Q$.
\end{proposition}

In the above LC circuit example, the Kirchhoff current law constraints $\Delta_{Q}$ in Eq.~\eqref{eq:KCL} induce the Dirac structure $D_{\Delta_{Q}}$.
In coordinates, we write an element in $T^{*}Q$ as $(q,p)$ with $p = (p_{\ell}, p_{c_{1}}, p_{c_{2}}, p_{c_{3}})$ and then by noting that $\Omega = dq \wedge dp$, we have
\begin{equation*}
  D_{\Delta_{Q}}(q,p) =
  \setdef{ (\dot{q}, \dot{p}, \alpha_{q}, \alpha_{p}) \in TT^{*}Q \oplus T^{*}T^{*}Q }
  {
    \dot{q} \in \Delta_{Q},
    \;
    \dot{q} = \alpha_{p},
    \;
    \dot{p} + \alpha_{q} \in \Delta_{Q}^{\circ}
  },
\end{equation*}
where $\Delta_{Q}^{\circ} \subset T^{*}Q$ is the annihilator of $\Delta_{Q}$ defined in Eq.~\eqref{eq:Delta_Q^circ}.

\subsection{Tulczyjew's Triple}
\label{ssec:TulczyjewTriple}
Following \citet{Tu1976a, Tu1976b} and \citet{YoMa2006a}, let us introduce Tulczyjew's triple, i.e., the diffeomorphisms $\Omega^{\flat}$, $\kappa_{Q}$, and $\gamma_{Q} \defeq \Omega^{\flat} \circ \kappa_{Q}^{-1}$ defined between the iterated tangent and cotangent bundles as follows.
\begin{subequations}
  \label{eq:TulczyjewTriple}
  \begin{equation}
    \vcenter{
      \xymatrix@!0@R=0.75in@C=.75in{
        T^{*}TQ \ar[dr]_{\pi_{TQ}\!\!} \ar@/^{1.75pc}/[rrrr]^{\gamma_{Q}} & &
        TT^{*}Q \ar[rr]^{\Omega^{\flat}} \ar[ll]_{\kappa_{Q}} \ar[dr]_{\tau_{T^{*}Q}\!\!} \ar[dl]^{\!\!T\pi_{Q}} & &
        T^{*}T^{*}Q \ar[dl]^{\!\!\pi_{T^{*}Q}}
        \\
         & TQ & & T^{*}Q &
      }
    }
  \end{equation}
  \begin{equation}
    \vcenter{
      \xymatrix@!0@R=0.75in@C=.75in{
        (q, \delta q, \delta p, p) \ar@{|->}[dr] & & 
        (q, p, \delta q, \delta p) \ar@{|->}[ll] \ar@{|->}[rr] \ar@{|->}[dl] \ar@{|->}[dr] & &
        (q, p, -\delta p, \delta q) \ar@{|->}[dl]
        \\
         & (q, \delta q) & & (q,p) &
      }
    }
  \end{equation}
\end{subequations}
The maps $\Omega^{\flat}$ and $\kappa_{Q}$ induce symplectic forms on $TT^{*}Q$ in the following way:
Let $\Theta_{T^{*}T^{*}Q}$ and $\Theta_{T^{*}TQ}$ be standard symplectic one-forms on the cotangent bundles $T^{*}T^{*}Q$ and $T^{*}TQ$, respectively.
One defines one-forms $\chi$ and $\lambda$ on $TT^{*}Q$ by
\begin{equation*}
  \chi \defeq (\Omega^{\flat})^{*} \Theta_{T^{*}T^{*}Q} = -\delta p\,dq + \delta q\,dp,
  \qquad
  \lambda \defeq (\kappa_{Q})^{*} \Theta_{T^{*}TQ} = \delta p\,dq + p\,d(\delta q),
\end{equation*}
and, using these one-forms, define the two-from $\Omega_{TT^{*}Q}$ on $TT^{*}Q$ by
\begin{equation*}
  \Omega_{TT^{*}Q} \defeq -d\lambda = d\chi = dq \wedge d(\delta p) + d(\delta q) \wedge dp.
\end{equation*}
Then, this gives a symplectic form on $TT^{*}Q$.

\subsection{Lagrange--Dirac Systems}
To define a Lagrange--Dirac system, it is necessary to introduce the Dirac differential of a Lagrangian function:
Given a Lagrangian $L:TQ \to \R$, we define the {\em Dirac differential} $\mathfrak{D}L: TQ \to T^*T^*Q$ by
\begin{equation*}
  \mathfrak{D}L \defeq \gamma_{Q} \circ dL.
\end{equation*}
In local coordinates,
\begin{equation*}
  \mathfrak{D}L(q,v) = 
  \parentheses{
    q, \pd{L}{v}, -\pd{L}{q}, v
  }.
\end{equation*}
Now we are ready to define a Lagrange--Dirac system:
\begin{definition}[Lagrange--Dirac Systems]
  \label{def:LDS}
  Suppose that a Lagrangian $L: TQ \to \R$ and a Dirac structure $D \subset TT^{*}Q \oplus T^{*}T^{*}Q$ are given.
  Let $X \in \mathfrak{X}(T^{*}Q)$ be a vector field on $T^{*}Q$.
  Then a {\em Lagrange--Dirac system} is defined by
  \begin{equation}
    \label{eq:LDS}
    (X, \mathfrak{D}L) \in D.
  \end{equation}
\end{definition}
In particular, if $D$ is the induced Dirac structure $D_{\Delta_{Q}}$ given in Eq.~\eqref{eq:D_Delta_Q}, the Lagrange--Dirac system can be written as follows:
\begin{equation*}
  T\pi_{Q}(X) \in \Delta_{Q},
  \qquad
  \Omega^{\flat}(X) - \mathfrak{D}L \in \Delta_{T^{*}Q}^{\circ},
\end{equation*}
or in local coordinates, by setting $X = \dot{q}\,\partial_{q} + \dot{p}\,\partial_{p}$, 
\begin{equation}
  \label{eq:LDS-coord}
  \dot{q} \in \Delta_{Q},
  \qquad
  \dot{q} = v,
  \qquad
  p = \pd{L}{v},
  \qquad
  \dot{p} - \pd{L}{q} \in \Delta_{Q}^{\circ}.
\end{equation}

\begin{example}[LC circuit]
  With the Dirac structure $D_{\Delta_{Q}}$ in Eq.~\eqref{eq:D_Delta_Q} induced by the constraints $\Delta_{Q}$ in Eq.~\eqref{eq:KCL}, the Lagrange--Dirac system $(X, \mathfrak{D}L) \in D_{\Delta_{Q}}$ gives
  \begin{subequations}
    \label{eq:LDS-coord-LC}
    \begin{equation}
      \dot{q} \in \Delta_{Q},
      \qquad
      \dot{q} = f,
      \qquad
      p = \pd{L}{f},
      \qquad
      \dot{p} - \pd{L}{q} = \mu_{1} \omega^{1} + \mu_{2} \omega^{2}
    \end{equation}
    with the Lagrange multipliers $\mu_{1}, \mu_{2} \in \R$; to be
    more explicit,
    \begin{equation}
      \begin{array}{c}
        \dot{q}^{\ell} = \dot{q}^{c_{2}},
        \quad
        \dot{q}^{c_{1}} = \dot{q}^{c_{2}} - \dot{q}^{c_{3}},
        \medskip\\
        \dot{q}^{\ell} = f^{\ell},
        \quad
        \dot{q}^{c_{1}} = f^{c_{1}},
        \quad
        \dot{q}^{c_{2}} = f^{c_{2}},
        \quad
        \dot{q}^{c_{3}} = f^{c_{3}},
        \medskip\\
        p_{\ell} = \ell\,f^{\ell},
        \quad
        p_{c_{1}} = p_{c_{2}} = p_{c_{3}} = 0,
        \medskip\\
        \displaystyle
        \dot{p}_{\ell} = -\mu_{1},
        \quad
        \dot{p}_{c_{1}} + \frac{q^{c_{1}}}{c_{1}} = -\mu_{2},
        \quad
        \dot{p}_{c_{2}} + \frac{q^{c_{2}}}{c_{2}} = \mu_{1} + \mu_{2},
        \quad
        \dot{p}_{c_{3}} + \frac{q^{c_{3}}}{c_{3}} = -\mu_{2}.
      \end{array}
    \end{equation}
  \end{subequations}
  This formulation recovers the equations given by circuit theory.
\end{example}

\begin{remark}
  Notice that this formulation by \citet{YoMa2006a} does not use the Kirchhoff voltage law; it instead uses the Kirchhoff current law with the symplectic structure on $T^{*}Q$ to define the Dirac structure $D_{\Delta_{Q}} \subset TT^{*}Q \oplus T^{*}T^{*}Q$.
  On the other hand, the formulation by \citet{BlCr1997} and \citet{Sc1998} uses the Dirac structure $D \subset TP \oplus T^{*}P$, with a different configuration space $P$, defined by {\em both} the Kirchhoff voltage and current laws, without using any additional geometric (symplectic) structure.
\end{remark}

\subsection{Implicit and Nonholonomic Hamiltonian Systems}
One can define an implicit Hamiltonian system in an analogous way as shown by \citet{Sc1998} and \citet{DaSc1998}:
\begin{definition}
  Suppose that a Hamiltonian $H: T^{*}Q \to \R$ and a Dirac structure $D \subset TT^{*}Q \oplus T^{*}T^{*}Q$ are given.
  Let $X \in \mathfrak{X}(T^{*}Q)$ be a vector field on $T^{*}Q$.
  Then an {\em implicit Hamiltonian system (IHS)} is defined by
  \begin{equation*}
    (X, dH) \in D.
  \end{equation*}
\end{definition}
In particular, if $D$ is the induced Dirac structure $D_{\Delta_{Q}}$ given in Eq.~\eqref{eq:D_Delta_Q}, the IHS gives the {\em nonholonomic Hamilton's equations}~(see, e.g., \citet{BaSn1993}, \citet{ScMa1994}, and \citet{KoMa1997c}):
\begin{equation*}
  T\pi_{Q}(X) \in \Delta_{Q},
  \qquad
  \Omega^{\flat}(X) - dH \in \Delta_{T^{*}Q}^{\circ},
\end{equation*}
or in local coordinates, by setting $X = \dot{q}\,\partial_{q} + \dot{p}\,\partial_{p}$, 
\begin{equation*}
  \dot{q} \in \Delta_{Q},
  \qquad
  \dot{q} = \pd{H}{p},
  \qquad
  \dot{p} + \pd{H}{q} \in \Delta_{Q}^{\circ}.
\end{equation*}
To keep the exposition in this section concise, we will not go into details about IHS here.
We would like to point the reader to the references cited above for details and examples of IHS.

\section{Discrete Analogues of Tulczyjew's Triple}
\label{sec:DATT}
In this section, we construct discrete analogues of Tulczyjew's triple shown in Eq.~\eqref{eq:TulczyjewTriple} that retain the key geometric properties, especially the symplecticity of the maps involved.
This makes it possible to formulate a natural structure-preserving discrete analogue of Lagrange--Dirac systems.
The discussion here is limited to the case where the configuration space $Q$ is a vector space.

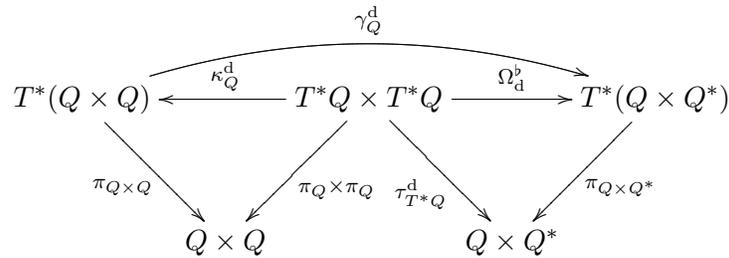
\begin{figure}[ht!]
  \centering
  $\xymatrix@!0@R=0.75in@C=.75in{
    T^{*}(Q \times Q) \ar[dr]_{\pi_{Q \times Q}\!\!} \ar@/^{1.75pc}/[rrrr]^{\gamma_{Q}^{\rm d}} & &
    T^{*}Q \times T^{*}Q \ar[rr]^{\Omega^{\flat}_{\rm d}} \ar[ll]_{\kappa_{Q}^{\rm d}} \ar[dr]_{\tau_{T^{*}Q}^{\rm d}\!\!\!\!} \ar[dl]^{\!\!\pi_{Q}\times\pi_{Q}} & &
    T^{*}(Q \times Q^{*}) \ar[dl]^{\!\!\pi_{Q \times Q^{*}}}
    \\
    & Q \times Q & & Q \times Q^{*} &
  }$
  \caption{A rough picture of a discrete analogue of Tulczyjew's triple.}
  \label{fig:DiscreteTulczyjewTriple}
\end{figure}

To give the big picture of what we would like to do in this section, constructing a discrete analogue of Tulczyjew's triple involves replacing, for example, the tangent bundle $TQ$ in Eq.~\eqref{eq:TulczyjewTriple} by the product $Q \times Q$  in accordance with the basic idea of discrete mechanics~(see, e.g., \cite{MaWe2001}); likewise $TT^{*}Q$ is replaced by $T^{*}Q \times T^{*}Q$; the role of $T^{*}Q$ in discrete mechanics is quite subtle in general, but since $Q$ is assumed to be a vector space, we can replace it with $Q \times Q^{*}$.
Fig.~\ref{fig:DiscreteTulczyjewTriple} gives a rough picture of a discrete analogue of Tulczyjew's triple. We work out the details of how to obtain the maps $\kappa_{Q}^{\rm d}$ and $\Omega^{\flat}_{\rm d}$ in the sections to follow.
The guiding principle here is to make use of symplectic maps associated with generating functions instead of smooth symplectic flows.

\subsection{Discrete Mechanics and Generating Functions}
\label{ssec:DiscreteMechAndGenFunctions}
Let us first review some basic facts on generating functions.
Consider a map $F: T^{*}Q \to T^{*}Q$ written as $(q_{0}, p_{0}) \mapsto (q_{1}, p_{1})$.
Note that, since $Q$ is assumed to be a vector space here, the cotangent bundle is trivial, i.e., $T^{*}Q \cong Q \times Q^{*}$, and so one can write $F: Q \times Q^{*} \to Q \times Q^{*}$ as well.
One then considers the following four maps associated with $F$:
\begin{enumerate}
\renewcommand{\theenumi}{\roman{enumi}}
\renewcommand{\labelenumi}{(\theenumi)}
\item $F_{1}: Q \times Q \to Q^{*} \times Q^{*};\ (q_{0}, q_{1}) \mapsto (p_{0}, p_{1})$,
\item $F_{2}: Q \times Q^{*} \to Q^{*} \times Q;\ (q_{0}, p_{1}) \mapsto (p_{0}, q_{1})$,
\item $F_{3}: Q^{*} \times Q \to Q \times Q^{*};\ (p_{0}, q_{1}) \mapsto (q_{0}, p_{1})$,
\item $F_{4}: Q^{*} \times Q^{*} \to Q \times Q;\ (p_{0}, p_{1}) \mapsto (q_{0}, q_{1})$.
\end{enumerate}
The Type~$i$ generating function with $i = 1, 2, 3, 4$ (using the terminology set by \citet{GoPoSa2001}) is a scalar function $S_{i}$ defined on the range of the map $F_{i}$ that exists if and only if the map $F$ is symplectic.
Let us look at the first three cases (the fourth one is not important here) and their relationship to discrete analogues of the map $\kappa_{Q}$ and $\Omega^{\flat}$ in the sections to follow.

\subsection{Generating Function of Type 1 and the Map $\kappa_{Q}^{\rm d}$}
This section relates the Type~1 generating function with a discrete analogue $\kappa_{Q}^{\rm d}$ of the map $\kappa_{Q}$ in Tulczyjew's triple, Eq.~\eqref{eq:TulczyjewTriple}.

First, we regard $(p_{0}, p_{1})$ as functions of $(q_{0}, q_{1})$ as indicated in the definition of the map $F_{1}$ above, and then define $i_{F_{1}}: Q \times Q \to T^{*}Q \times T^{*}Q$ by
\begin{equation*}
  i_{F_{1}}: (q_{0}, q_{1}) \mapsto ( (q_{0}, p_{0}), (q_{1}, p_{1}) )
  \quad
  \text{where}
  \quad
  (p_{0}, p_{1}) = F_{1}(q_{0}, q_{1}).
\end{equation*}
Now recall that the map $F: (q_{0}, p_{0}) \mapsto (q_{1}, p_{1})$ is symplectic if and only if $dq_{0} \wedge dp_{0} = dq_{1} \wedge dp_{1}$, or equivalently $d( -p_{0}\,dq_{0} + p_{1}\,dq_{1} ) = 0$.
Then, the Poincar\'e lemma states that this is true if and only if there exists some function $S_{1}: Q \times Q \to \R$, a {\em Type~1 generating function}, such that
\begin{equation*}
  -p_{0}\,dq_{0} + p_{1}\,dq_{1} = dS_{1}(q_{0}, q_{1}).
\end{equation*}
This relates the $(p_{0}, p_{1})$ with the generating function $S_{1}$:
\begin{equation}
  \label{eq:S_{1}eq-vecsp}
  p_{0} = -D_{1}S_{1}(q_{0}, q_{1}),
  \qquad
  p_{1} = D_{2}S_{1}(q_{0}, q_{1}).
\end{equation}
Then, this gives rise to the map $\kappa_{Q}^{\rm d}: T^{*}Q \times T^{*}Q \to T^{*}(Q \times Q)$ so that the diagram
\begin{equation}
  \label{cd:kappa_Q^d-vecsp}
  \begin{array}{cc}
    \xymatrix@!0@R=.7in@C=.8in{
      T^{*}Q \times T^{*}Q \ar[rr]^{\kappa_{Q}^{\rm d}} & & T^{*}(Q \times Q)
      \\
      & Q \times Q \ar[ul]^{i_{F_{1}}} \ar[ur]_{dS_{1}} & 
    }
    &
    \xymatrix@!0@R=.7in@C=.8in{
      ((q_{0},p_{0}),(q_{1},p_{1})) \ar@{|->}[rr] & & (q_{0}, q_{1}, -p_{0}, p_{1})
      \\
      & (q_{0},q_{1}) \ar@{|->}[ul] \ar@{|->}[ur] & 
    }
  \end{array}
\end{equation}
commutes, i.e., we obtain
\begin{equation}
  \label{eq:kappa_Q^d-vecsp}
  \kappa_{Q}^{\rm d}: ((q_{0},p_{0}),(q_{1},p_{1})) \mapsto (q_{0}, q_{1}, -p_{0}, p_{1}).
\end{equation}

\subsection{Generating Function of Type 2 and the Map $\Omega^{\flat}_{\rm d+}$}
Next, we would like to relate the Type~2 generating function with one of the two discrete analogues of the map $\Omega^{\flat}$ in Tulczyjew's triple, Eq.~\eqref{eq:TulczyjewTriple}.

First, we regard $(p_{0}, q_{1})$ as functions of $(q_{0}, p_{1})$ as indicated in the definition of the map $F_{2}$ above, and then define $i_{F_{2}}: Q \times Q^{*} \to T^{*}Q \times T^{*}Q$ by
\begin{equation*}
  i_{F_{2}}: (q_{0}, p_{1}) \mapsto ( (q_{0}, p_{0}), (q_{1}, p_{1}) )
  \quad
  \text{where}
  \quad
  (p_{0}, q_{1}) = F_{2}(q_{0}, p_{1}).
\end{equation*}
The map $F: (q_{0}, p_{0}) \mapsto (q_{1}, p_{1})$ is symplectic if and only if $dq_{0} \wedge dp_{0} = dq_{1} \wedge dp_{1}$, or equivalently $d(p_{0}\,dq_{0} + q_{1}\,dp_{1}) = 0$.
Then, the Poincar\'e lemma states that this is true if and only if there exists some function $S_{2}: Q \times Q^{*} \to \R$, a {\em Type~2 generating function}, such that
\begin{equation*}
  p_{0}\,dq_{0} + q_{1}\,dp_{1} = dS_{2}(q_{0}, p_{1}).
\end{equation*}
This relates the $(p_{0}, q_{1})$ with the generating function $S_{2}$:
\begin{equation}
  \label{eq:S_{2}eq-vecsp}
  p_{0} = D_{1}S_{2}(q_{0}, p_{1}),
  \qquad
  q_{1} = D_{2}S_{2}(q_{0}, p_{1}).
\end{equation}
Then, this gives rise to the map $\Omega^{\flat}_{\rm d+}: T^{*}Q \times T^{*}Q \to T^{*}(Q \times Q^{*})$ so that the diagram
\begin{equation}
  \label{cd:Omega_d+^flat-vecsp}
  \begin{array}{cc}
    \xymatrix@!0@R=0.7in@C=.8in{
      T^{*}Q \times T^{*}Q \ar[rr]^{\Omega^{\flat}_{\rm d+}} & & T^{*}(Q \times Q^{*})
      \\
      & Q \times Q^{*} \ar[ul]^{i_{F_{2}}} \ar[ur]_{dS_{2}} & 
    }
    &
    \xymatrix@!0@R=0.7in@C=.8in{
      ((q_{0},p_{0}),(q_{1},p_{1})) \ar@{|->}[rr] & & (q_{0}, p_{1}, p_{0}, q_{1})
      \\
      & (q_{0},p_{1}) \ar@{|->}[ul] \ar@{|->}[ur] & 
    }
  \end{array}
\end{equation}
commutes, i.e., we obtain
\begin{equation}
  \label{eq:Omega_d+^flat-vecsp}
  \Omega^{\flat}_{\rm d+}: ((q_{0},p_{0}),(q_{1},p_{1})) \mapsto (q_{0}, p_{1}, p_{0}, q_{1}).
\end{equation}

\subsection{Generating Function of Type 3 and the Map $\Omega^{\flat}_{\rm d-}$}
The other discrete analogue of the map $\Omega^{\flat}$ follows from the Type~3 generating function.

In this case, we regard $(q_{0}, p_{1})$ as functions of $(p_{0}, q_{1})$ as indicated in the definition of the map $F_{3}$ above, and then define $i_{F_{3}}: Q^{*} \times Q \to T^{*}Q \times T^{*}Q$ by
\begin{equation*}
  i_{F_{3}}: (p_{0}, q_{1}) \mapsto ( (q_{0}, p_{0}), (q_{1}, p_{1}) )
  \quad
  \text{where}
  \quad
  (q_{0}, p_{1}) = F_{3}(p_{0}, q_{1}).
\end{equation*}
The map $F: (q_{0}, p_{0}) \mapsto (q_{1}, p_{1})$ is symplectic if and only if $dq_{0} \wedge dp_{0} = dq_{1} \wedge dp_{1}$, or equivalently $d(-q_{0}\,dp_{0} - p_{1}\,dq_{1}) = 0$.
Then, again by the Poincar\'e lemma, this is true if and only if there exists some function $S_{3}: Q^{*} \times Q \to \R$ such that
\begin{equation*}
  -q_{0}\,dp_{0} - p_{1}\,dq_{1} = dS_{3}(p_{0}, q_{1}).
\end{equation*}
This relates the $(q_{0}, p_{1})$ with the generating function $S_{3}$:
\begin{equation}
  \label{eq:S_{3}eq-vecsp}
  q_{0} = -D_{1}S_{3}(p_{0}, q_{1}),
  \qquad
  p_{1} = -D_{2}S_{3}(p_{0}, q_{1}).
\end{equation}
Then, this gives rise to the map $\Omega^{\flat}_{\rm d-}: T^{*}Q \times T^{*}Q \to T^{*}(Q^{*} \times Q)$ so that the diagram
\begin{equation}
  \label{cd:Omega_d-^flat-vecsp}
  \begin{array}{cc}
    \xymatrix@!0@R=0.7in@C=.8in{
      T^{*}Q \times T^{*}Q \ar[rr]^{\Omega^{\flat}_{\rm d-}} & & T^{*}(Q^{*} \times Q)
      \\
      & Q^{*} \times Q \ar[ul]^{i_{F_{3}}} \ar[ur]_{dS_{3}} & 
    }
    &
    \xymatrix@!0@R=0.7in@C=.8in{
      ((q_{0},p_{0}),(q_{1},p_{1})) \ar@{|->}[rr] & & (p_{0}, q_{1}, -q_{0}, -p_{1})
      \\
      & (p_{0},q_{1}) \ar@{|->}[ul] \ar@{|->}[ur] & 
    }
  \end{array}
\end{equation}
commutes, i.e., we obtain
\begin{equation}
  \label{eq:Omega_d-^flat-vecsp}
  \Omega^{\flat}_{\rm d-}: ((q_{0},p_{0}),(q_{1},p_{1})) \mapsto (p_{0}, q_{1}, -q_{0}, -p_{1}).
\end{equation}

\subsection{$(+)$-Discrete Tulczyjew Triple}
\label{ssec:+DTT-vecsp}
Combining the diagrams in Eqs.~\eqref{cd:kappa_Q^d-vecsp} and \eqref{cd:Omega_d+^flat-vecsp}, we obtain the following {\em $(+)$-discrete Tulczyjew triple}.
\begin{subequations}
  \begin{equation}
    \label{eq:+DTulczyjewTriple-vecsp}
    \vcenter{
      \xymatrix@!0@R=0.75in@C=.75in{
        T^{*}(Q \times Q) \ar[dr]_{\pi_{Q \times Q}\!\!} \ar@/^{1.75pc}/[rrrr]^{\gamma_{Q}^{\rm d+}} & &
        T^{*}Q \times T^{*}Q \ar[rr]^{\Omega^{\flat}_{\rm d+}} \ar[ll]_{\kappa_{Q}^{\rm d}} \ar[dr]_{\tau_{T^{*}Q}^{\rm d+}\!\!\!\!} \ar[dl]^{\!\!\pi_{Q} \times \pi_{Q}} & &
        T^{*}(Q \times Q^{*}) \ar[dl]^{\!\!\pi_{Q \times Q^{*}}}
        \\
        & Q \times Q & & Q \times Q^{*} &
      }
    }
  \end{equation}
  \begin{equation}
    \vcenter{
      \xymatrix@!0@R=0.75in@C=.75in{
        (q_{0}, q_{1}, -p_{0}, p_{1}) \ar@{|->}[dr] & & 
        ((q_{0}, p_{0}), (q_{1}, p_{1})) \ar@{|->}[ll] \ar@{|->}[rr] \ar@{|->}[dl] \ar@{|->}[dr] & &
        (q_{0}, p_{1}, p_{0}, q_{1}) \ar@{|->}[dl]
        \\
        & (q_{0}, q_{1}) & & (q_{0}, p_{1}) &
      }
    }
  \end{equation}
\end{subequations}

The maps $\kappa_{Q}^{\rm d}$ and $\Omega^{\flat}_{\rm d+}$ inherit the properties of $\kappa_{Q}$ and $\Omega^{\flat}$ discussed in Section~\ref{ssec:TulczyjewTriple} in the following sense:
Let $\Theta_{T^{*}(Q \times Q^{*})}$ and $\Theta_{T^{*}(Q \times Q)}$ be the symplectic one-forms on $T^{*}(Q \times Q^{*})$ and $T^{*}(Q \times Q)$, respectively.
The maps $\kappa_{Q}^{\rm d}$ and $\Omega^{\flat}_{\rm d+}$ induce two symplectic one-forms on $T^{*}Q \times T^{*}Q$.
One is
\begin{equation*}
  \chi_{\rm d+} \defeq (\Omega^{\flat}_{\rm d+})^{*} \Theta_{T^{*}(Q \times Q^{*})} = p_{0}\,dq_{0} + q_{1}\,dp_{1},
\end{equation*}
and the other is
\begin{equation*}
  \lambda_{\rm d+} \defeq (\kappa_{Q}^{\rm d})^{*} \Theta_{T^{*}(Q \times Q)} = -p_{0}\,dq_{0} + p_{1}\,dq_{1}.
\end{equation*}
Then, using these one-forms, define the two-from $\Omega_{T^{*}Q \times T^{*}Q}$ by
\begin{equation*}
  \Omega_{T^{*}Q \times T^{*}Q} = -d\lambda_{\rm d+} = d\chi_{\rm d+} = dq_{1} \wedge dp_{1} - dq_{0} \wedge dp_{0}.
\end{equation*}
This is a natural symplectic form defined on the product of two cotangent bundles (see \citet[Proposition~5.2.1 on p.~379]{AbMa1978}).

\subsection{$(-)$-Discrete Tulczyjew Triple}
\label{ssec:-DTT-vecsp}
Combining the diagrams in Eqs.~\eqref{cd:kappa_Q^d-vecsp} and \eqref{cd:Omega_d-^flat-vecsp}, we obtain the following {\em $(-)$-discrete Tulczyjew triple}.
\begin{subequations}
  \begin{equation}
    \label{eq:-DTulczyjewTriple-vecsp}
    \vcenter{
      \xymatrix@!0@R=0.7in@C=.7in{
        T^{*}(Q \times Q) \ar[dr]_{\pi_{Q \times Q}\!\!} \ar@/^{1.75pc}/[rrrr]^{\gamma_{Q}^{\rm d-}} & &
        T^{*}Q \times T^{*}Q \ar[rr]^{\Omega^{\flat}_{\rm d-}} \ar[ll]_{\kappa_{Q}^{\rm d}} \ar[dr]_{\tau_{T^{*}Q}^{\rm d-}\!\!\!\!} \ar[dl]^{\!\!\pi_{Q} \times \pi_{Q}} & &
        T^{*}(Q^{*} \times Q) \ar[dl]^{\!\!\pi_{Q^{*} \times Q}}
        \\
        & Q \times Q & & Q^{*} \times Q &
      }
    }
  \end{equation}
  \begin{equation}
    \vcenter{
      \xymatrix@!0@R=0.7in@C=.7in{
        (q_{0}, q_{1}, -p_{0}, p_{1}) \ar@{|->}[dr] & & 
        ((q_{0}, p_{0}), (q_{1}, p_{1})) \ar@{|->}[ll] \ar@{|->}[rr] \ar@{|->}[dl] \ar@{|->}[dr] & &
        (p_{0}, q_{1}, -q_{0}, -p_{1}) \ar@{|->}[dl]
        \\
        & (q_{0}, q_{1}) & & (p_{0}, q_{1}) &
      }
    }
  \end{equation}
\end{subequations}

As in the $(+)$-discrete case, the maps $\kappa_{Q}^{\rm d}$ and $\Omega^{\flat}_{\rm d-}$ inherit the properties of $\kappa_{Q}$ and $\Omega^{\flat}$:
Let $\Theta_{T^{*}(Q^{*} \times Q)}$ be the symplectic one-form on $T^{*}(Q^{*} \times Q)$.
Then, we have
\begin{equation*}
  \chi_{\rm d-} \defeq (\Omega^{\flat}_{\rm d-})^{*} \Theta_{T^{*}(Q^{*} \times Q)} = -p_{1}\,dq_{1} - q_{0}\,dp_{0},
\end{equation*}
and
\begin{equation*}
  \lambda_{\rm d-} \defeq (\kappa_{Q}^{\rm d})^{*} \Theta_{T^{*}(Q \times Q)} = -p_{0}\,dq_{0} + p_{1}\,dq_{1}.
\end{equation*}
Then, they induce the same symplectic form $\Omega_{T^{*}Q \times T^{*}Q}$ as above:
\begin{equation*}
  \Omega_{T^{*}Q \times T^{*}Q} \defeq -d\lambda_{\rm d-} = d\chi_{\rm d-} = dq_{1} \wedge dp_{1} - dq_{0} \wedge dp_{0}.
\end{equation*}

\section{Discrete Analogues of Induced Dirac Structures}
\label{sec:DAIDS}
Recall from Section~\ref{ssec:Induced_Dirac_Structure} that, given a constraint distribution $\Delta_{Q} \subset TQ$, we first defined the distribution $\Delta_{T^{*}Q} \subset TT^{*}Q$ and then constructed the induced Dirac structure $D_{\Delta_{Q}} \subset TT^{*}Q \oplus T^{*}T^{*}Q$.
This section develops a discrete analogue of this construction.
\subsection{Discrete constraint distributions}
\label{subsec:discrete_constraints}
Given the fact that the tangent bundle $TQ$ is replaced by the product $Q \times Q$ in the discrete setting, a natural discrete analogue of a constraint distribution $\Delta_{Q} \subset TQ$ is a subset $\Delta_{Q}^{\rm d} \subset Q \times Q$.
We follow the approach of \citet{CoMa2001} (see also \citet{McPe2006}) to construct discrete constraints $\Delta_{Q}^{\rm d} \subset Q \times Q$ based on given (continuous) constraints $\Delta_{Q} \subset TQ$.

Let $\Delta_{Q}^{\circ} \subset T^{*}Q$ be the annihilator distribution (or codistribution) of $\Delta_{Q} \subset TQ$ and $m \defeq \dim T_{q}Q - \dim\Delta_{Q}(q)$ for each $q \in Q$.
Then, one can find a set of $m$ constraint one-forms $\{ \omega^{a} \}_{a=1}^{m}$ that spans the annihilator:
\begin{equation*}
  \Delta_{Q}^{\circ} = \Span\{ \omega^{a} \}_{a=1}^{m}.
\end{equation*}
In local coordinates, we may write
\begin{equation}
  \label{eq:omega-coord}
  \omega^{a}(q, v) = A_{i}^{a}(q) v^{i},
\end{equation}
where $(A^{a}_{i}(q))$ is an $m \times n$ full-rank matrix for each $q \in Q$, i.e., $\rank A(q) = m$.

Then, using the one-forms $\omega^{a}$ and a retraction $\mathcal{R}: TQ \to Q$ (see Section~\ref{ssec:Retractions}), we define functions $\omega_{\rm d\pm}^{a}: Q \times Q \to \R$ by
\begin{equation}
  \label{eq:omega_d+-}
  \omega_{\rm d+}^{a}(q_{0}, q_{1}) \defeq \omega^{a}\parentheses{ q_{0}, \mathcal{R}_{q_{0}}^{-1}(q_{1}) },
  \qquad
  \omega_{\rm d-}^{a}(q_{0}, q_{1}) \defeq \omega^{a}\parentheses{ q_{1}, -\mathcal{R}_{q_{1}}^{-1}(q_{0}) },
\end{equation}
and then define the discrete constraints $\Delta_{Q}^{\rm d\pm} \subset Q \times Q$ as follows:
\begin{equation}
  \label{eq:Delta_Q^d}
  \Delta_{Q}^{\rm d\pm} \defeq \setdef{ (q_{0}, q_{1}) \in Q \times Q }{ \omega_{\rm d\pm}^{a}(q_{0}, q_{1}) = 0,\, a = 1, 2, \dots, m }.
\end{equation}
The following proposition suggests that it is natural to think of $q_{1}$ as a discrete analogue of the velocity $v_{q_{0}} \in T_{q_{0}}Q$ when imposing the constraint $(q_{0}, q_{1}) \in \Delta_{Q}^{\rm d+}$, and $q_{0}$ a discrete analogue of $v_{q_{1}} \in T_{q_{1}}Q$ when imposing $(q_{0}, q_{1}) \in \Delta_{Q}^{\rm d-}$:
\begin{proposition}
  \label{prop:Delta_Q^d}
  The discrete constraints defined by $(q_{0}, q_{1}) \in \Delta_{Q}^{\rm d\pm} \subset Q \times Q$ are constraints only on the variable $q_{1}$ and $q_{0}$, respectively; i.e., $pr_{1}( \Delta_{Q}^{\rm d+} ) = Q$ and $pr_{2}( \Delta_{Q}^{\rm d-} ) = Q$, where $pr_{i}: Q \times Q \to Q$ with $i = 1, 2$ is the projection to the $i$-th component.
\end{proposition}

\begin{proof}
  Let $\omega: TQ \to \R^{m}$ be the map defined by 
  \begin{equation*}
    \omega(q, v) \defeq
    \parentheses{
      \omega^{1}(q, v), \dots, \omega^{m}(q, v)
    },
  \end{equation*}
  and $\omega_{\rm d+}: Q \times Q \to \R^{m}$ be the map defined by 
  \begin{equation*}
    \omega_{\rm d+}(q_{0}, q_{1}) \defeq
    \parentheses{
      \omega_{\rm d+}^{1}(q_{0}, q_{1}), \dots, \omega_{\rm d+}^{m}(q_{0}, q_{1})
    }.
  \end{equation*}
  In the first equation in Eq.~\eqref{eq:omega_d+-}, take the derivative respect to $q_{1}$ to obtain
  \begin{align*}
    D_{2}\omega_{\rm d+}(q_{0}, q_{1})
    &= D_{2}\omega\parentheses{ q_{0}, \mathcal{R}_{q_{0}}^{-1}(q_{1}) } \cdot D\mathcal{R}_{q_{0}}^{-1}(q_{1})
    \\
    &= A(q_{0}) \cdot D\mathcal{R}_{q_{0}}^{-1}(q_{1}),
  \end{align*}
  where we used the coordinate expression for $\omega$ in Eq.~\eqref{eq:omega-coord}.
  Since $D\mathcal{R}_{q_{0}}^{-1}$ is an invertible matrix (see Remark~\ref{rem:mathcalR_inverse}) and $\rank A = m$, we find that $\rank D_{2}\omega_{\rm d+} = m$.
  Therefore, by the implicit function theorem, we may (locally) rewrite the constraints $\omega_{\rm d+}(q_{0}, q_{1}) = 0$ as $q_{1}^{i_{l}} = f^{l}(q_{0}, q_{1}^{j_{1}}, \dots q_{1}^{j_{n-m}})$ with some function $f^{l}: \R^{n} \times \R^{n-m} \to \R^{m}$ for $l = 1, \dots, m$, where $\{ i_{1}, \dots, i_{m} \} \cup \{j_{1}, \dots, j_{n-m}\} = \{ 1, 2, \dots, n \}$ and $\{ i_{1}, \dots, i_{m} \} \cap \{j_{1}, \dots, j_{n-m}\} = \varnothing$.
  Hence $q_{0}$ is a free variable and so the claim follows.
  Similarly for $\omega_{\rm d-}$.
\end{proof}

Next, we introduce discrete analogues of the distribution $\Delta_{T^{*}Q} \subset TT^{*}Q$ using the discrete constraint $\Delta_{Q}^{\rm d\pm}$ defined above.
Natural discrete analogues of $\Delta_{T^{*}Q}$ would be $\Delta_{T^*Q}^{\rm d\pm} \subset T^{*}Q \times T^{*}Q$ defined by
\begin{equation*}
  \Delta_{T^*Q}^{\rm d\pm}
  \defeq (\pi_{Q}\times \pi_{Q})^{-1}(\Delta_{Q}^{\rm d\pm})
  = \setdef{\parentheses{(q_0, p_0),(q_1,p_1)}\in T^*Q\times T^*Q}{(q_0, q_1)\in\Delta_{Q}^{\rm d\pm}},
\end{equation*}
which is analogous to the continuous distribution $\Delta_{T^{*}Q} \defeq (T\pi_{Q})^{-1}(\Delta_{Q})$ in Eq.~\eqref{eq:Delta_T^*Q}.

We will also need discrete analogues of the annihilator $\Delta_{T^{*}Q}^{\circ}$ defined in Eq.~\eqref{eq:Delta_T^*Q}; natural discrete analogues of it would be annihilator distributions on $Q \times Q^{*}$ and $Q^{*} \times Q$.
We use the projections $\pi_{Q}^{\rm d+}: Q \times Q^{*} \to Q$ and $\pi_{Q}^{\rm d-}: Q^{*} \times Q \to Q$ to define annihilator distributions $\Delta_{Q \times Q^{*}}^{\circ} \subset T^{*}(Q \times Q^{*})$ and $\Delta_{Q^{*} \times Q}^{\circ} \subset T^{*}(Q^{*} \times Q)$ as follows:
\begin{align*}
  \Delta_{Q \times Q^{*}}^{\circ} & \defeq (\pi_{Q}^{\rm d+})^{*}(\Delta_{Q}^{\circ})
  = \setdef{(q, p, \alpha_{q}, 0) \in T^{*}(Q \times Q^{*})}{ \alpha_{q}\,dq \in \Delta_{Q}^{\circ}(q) },
  \\
  \Delta_{Q^{*} \times Q}^{\circ} & \defeq (\pi_{Q}^{\rm d-})^{*}(\Delta_{Q}^{\circ}) 
  =\setdef{(p, q, 0, \alpha_{q}) \in T^{*}(Q^{*} \times Q)}{ \alpha_{q}\,dq \in \Delta_{Q}^{\circ}(q) },
\end{align*}
which is analogous to the expression for the continuous annihilator distribution $\Delta_{T^{*}Q}^{\circ}=\pi_Q^*(\Delta_Q^\circ)$.

\subsection{Discrete Induced Dirac Structures}
Now we are ready to define discrete analogues of the induced Dirac structures $D_{\Delta_{Q}}$ shown in Proposition~\ref{prop:IDS}.
\begin{definition}[Discrete Induced Dirac Structures]
  Given a discrete constraint distribution $\Delta_{Q}^{\rm d+} \subset Q \times Q$, we define the {\em $(+)$-discrete induced Dirac structure} as follows:
  \begin{multline*}
    D_{\Delta_{Q}}^{\rm d+} \defeq
    \Bigl\{
    ((z,z^{+}), \alpha_{\hat{z}}) \in (T^{*}Q \times T^{*}Q) \times T^{*}(Q \times Q^{*})
    \  |\ \\
    \parentheses{z,z^{+}}\in\Delta_{T^*Q}^{\rm d+},\; \alpha_{\hat{z}} - \Omega^{\flat}_{\rm d+}\parentheses{z, z^{+}} \in \Delta_{Q \times Q^{*}}^\circ
  \Bigr\},
  \end{multline*}
  where if $z=(q,p)$ and $z^{+} = (q^{+}, p^{+})$ then $\hat{z} \defeq (q, p^{+}) \in Q \times Q^{*}$.
  Likewise, given a discrete constraint distribution $\Delta_{Q}^{\rm d-} \subset Q \times Q$, we define the {\em $(-)$-discrete induced Dirac structure} as follows:
  \begin{multline*}
    D_{\Delta_{Q}}^{\rm d-} \defeq
    \Bigl\{
      ((z^{-},z), \alpha_{\tilde{z}}) \in (T^{*}Q \times T^{*}Q) \times T^{*}(Q^{*} \times Q)
      \ |\ \\
      \parentheses{z^{-},z}\in\Delta_{T^*Q}^{\rm d-},\; \alpha_{\tilde{z}} - \Omega^{\flat}_{\rm d-}\parentheses{z^{-},z}\in \Delta_{Q^{*} \times Q}^\circ
    \Bigr\},
  \end{multline*}
  where if $z=(q, p)$ and $z^{-}=(q^{-}, p^{-})$ then $\tilde{z} \defeq (p^{-}, q) \in Q^{*} \times Q$.
\end{definition}

\section{Discrete Dirac Mechanics}
\label{sec:DDM}
Now that we have discrete analogues of both Tulczyjew's triple and induced Dirac structures at our disposal, we are ready to define discrete analogues of Lagrange--Dirac and nonholonomic Hamiltonian systems.
As we shall see, two types of discrete Lagrange--Dirac/nonholonomic Hamiltonian systems will follow from the $(\pm)$-discrete Tulczyjew triples and $(\pm)$-discrete induced Dirac structures.

\subsection{$(+)$-Discrete Dirac Mechanics}
\label{ssec:+DDS}
\subsubsection{$(+)$-Discrete Lagrange--Dirac Systems}
Let us first introduce a discrete analogue of the Dirac differential: Define $\gamma_{Q}^{\rm d+}: T^{*}(Q \times Q) \to T^{*}(Q \times Q^{*})$ by
\begin{equation*}
  \gamma_{Q}^{\rm d+} \defeq \Omega^{\flat}_{\rm d+} \circ (\kappa_{Q}^{\rm d})^{-1},
\end{equation*}
and, for a given discrete Lagrangian $L_{\rm d}: Q \times Q \to \R$, define the {\em $(+)$-discrete Dirac differential} $\mathfrak{D}^{+}L_{\rm d}: Q \times Q \to T^{*}(Q \times Q^{*})$ by
\begin{equation*}
  \mathfrak{D}^{+}L_{\rm d} \defeq \gamma_{Q}^{\rm d+} \circ dL_{\rm d}.
\end{equation*}
In coordinates, we have
\begin{equation*}
  \mathfrak{D}^{+}L_{\rm d}(q_{k}, q^{+}_{k}) = (q_{k}, D_{2}L_{\rm d}, -D_{1}L_{\rm d}, q^{+}_{k}).
\end{equation*}

\begin{definition}[$(+)$-Discrete Lagrange--Dirac System]
  \label{def:I+DLS}
  Suppose that a discrete Lagrangian $L_{\rm d}: Q \times Q \to \R$ and the constraint distribution $\Delta_{Q} \subset TQ$ are given; and so Eq.~\eqref{eq:Delta_Q^d} gives the discrete constraint distribution $\Delta_{Q}^{\rm d+} \subset Q \times Q$.
  Let
  \begin{equation}
    \label{eq:X_d}
    X_{\rm d}^{k} = ((q_{k}, p_{k}),(q_{k+1},p_{k+1})) \in T^{*}Q \times T^{*}Q
  \end{equation}
  be a discrete analogue of a vector field on $T^{*}Q$.
  Then, a {\em $(+)$-discrete Lagrange--Dirac system} is a triple $(L_{\rm d}, \Delta_{Q}, X_{\rm d})$ with
  \begin{equation}
    \label{eq:I+DLS}
    \parentheses{ X_{\rm d}^{k}, \mathfrak{D}^{+}L_{\rm d}(q_{k},q^{+}_{k}) } \in D_{\Delta_{Q}}^{\rm d+}.
  \end{equation}
\end{definition}

\begin{remark}
  The variable $q^{+}_{k}$ in Eq.~\eqref{eq:I+DLS} is a discrete analogue of $v$ in Eq.~\eqref{eq:LDS-coord}. See Proposition~\ref{prop:Delta_Q^d}.
\end{remark}

Let us find a coordinate expression for a $(+)$-discrete Lagrange--Dirac system: Eq.~\eqref{eq:I+DLS} gives
\begin{equation*}
  (q_k, q_{k+1}) \in \Delta_{Q}^{\rm d+},
  \qquad
  \mathfrak{D}^{+}L_{\rm d}- \Omega^{\flat}_{\rm d+}(X_{\rm d}^{k})\in \Delta_{Q \times Q^{*}}^\circ,
\end{equation*}
where
\begin{equation*}
  \Omega^{\flat}_{\rm d+}(X_{\rm d}^{k}) = (q_{k}, p_{k+1}, p_{k}, q_{k+1}).
\end{equation*}
Thus, we obtain the following set of equations:
\begin{subequations}
  \label{eq:+DLDEq}
  \begin{equation}
    \begin{array}{c}
      \DS (q_k,q_{k+1}) \in \Delta_{Q}^{\rm d+},
      \qquad
      \DS q_{k+1} = q^{+}_{k},
      \medskip\\
      \DS p_{k+1} = D_{2}L_{\rm d}(q_{k}, q^{+}_{k}),
      \qquad
      \DS p_{k} + D_{1}L_{\rm d}(q_{k}, q^{+}_{k}) \in \Delta_{Q}^{\circ}(q_{k}),
    \end{array}
  \end{equation}
  or more explicitly, with the Lagrange multipliers $\mu_{a}$,
  \begin{equation}
    \begin{array}{c}
      \DS \omega_{\rm d+}^{a}(q_k,q_{k+1}) = 0,
      \qquad
      \DS q_{k+1} = q^{+}_{k},
      \medskip\\
      \DS p_{k+1} = D_{2}L_{\rm d}(q_{k}, q^{+}_{k}),
      \qquad
      \DS p_{k} + D_{1}L_{\rm d}(q_{k}, q^{+}_{k}) = \mu_{a} \omega^{a}(q_{k}),
    \end{array}
  \end{equation}
\end{subequations}
where $a = 1, 2, \dots, m$.
We shall call them the {\em $(+)$-discrete Lagrange--Dirac equations}; they recover the nonholonomic integrator of \citet{CoMa2001} (see also \citet{McPe2006}).

Consider the special case $\Delta_{Q} = TQ$.
In this case, $\Delta_{Q}^{\rm d+} = Q \times Q$ and $\Delta_{Q}^{\circ} = 0$, and so the above equations reduce to
\begin{equation}
  \label{eq:DIELEq}
  q_{k+1} = q^{+}_{k},
  \qquad
  p_{k+1} = D_{2}L_{\rm d}(q_{k}, q^{+}_{k}),
  \qquad
  p_{k} = - D_{1}L_{\rm d}(q_{k}, q^{+}_{k}).
\end{equation}
These are equivalent to the discrete Euler--Lagrange equations (see \citet{MaWe2001}).

\subsubsection{$(+)$-Discrete Nonholonomic Hamiltonian System}
A nonholonomic discrete Hamiltonian system is defined analogously:
\begin{definition}[$(+)$-Discrete Nonholonomic Hamiltonian System]
  \label{def:I+DHS}
  Suppose that a $(+)$-discrete Hamiltonian~(referred to as the right discrete Hamiltonian in \citep{LaWe2006}) $H_{\rm d+}: Q \times Q^{*} \to \R$ and the constraint distribution $\Delta_{Q} \subset TQ$ are given; and so Eq.~\eqref{eq:Delta_Q^d} gives the discrete constraint distribution $\Delta_{Q}^{\rm d+} \subset Q \times Q$.
  Let $X_{\rm d}^{k}$ be a discrete analogue of a vector field on $T^{*}Q$ as in Eq.~\eqref{eq:X_d}.
  Then, a {\em $(+)$-discrete nonholonomic Hamiltonian system} is a triple $(H_{\rm d+}, \Delta_{Q}, X_{\rm d})$ with
  \begin{equation}
    \label{eq:I+DHS}
    \parentheses{X_{\rm d}^{k}, dH_{\rm d+}(q_{k}, p_{k+1})} \in D_{\Delta_{Q}}^{\rm d+}.
  \end{equation}
\end{definition}

A coordinate expression is obtained in a similar way:
\begin{subequations}
  \label{eq:I+DHS-coord}
  \begin{equation}
    (q_k, q_{k+1}) \in \Delta_{Q}^{\rm d+},
    \qquad
    q_{k+1} = D_{2}H_{\rm d+}(q_{k}, p_{k+1}),
    \qquad
    p_{k} - D_{1}H_{\rm d+}(q_{k}, p_{k+1}) \in \Delta_{Q}^{\circ}(q_{k}),
  \end{equation}
  or more explicitly
  \begin{equation}
    \omega_{\rm d+}^{a}(q_k,q_{k+1}) = 0,
    \qquad
    q_{k+1} = D_{2}H_{\rm d+}(q_{k}, p_{k+1}),
    \qquad
    p_{k} - D_{1}H_{\rm d+}(q_{k}, p_{k+1}) = \mu_{a} \omega^{a}(q_{k}),
  \end{equation}
\end{subequations}
where $a = 1, 2, \dots, m$.
We shall call them the {\em $(+)$-discrete nonholonomic Hamilton's equations}.

If $\Delta_{Q} = TQ$, then $\Delta_{Q}^{\rm d+} = Q\times Q$ and $\Delta_{Q}^{\circ} = 0$; and so the above equations reduce to
\begin{equation}
  q_{k+1} = D_{2}H_{\rm d+}(q_{k}, p_{k+1}),
  \qquad
  p_{k} = D_{1}H_{\rm d+}(q_{k}, p_{k+1}),
\end{equation}
which are the right discrete Hamilton's equations in Lall and West~\citep{LaWe2006}.

\subsection{$(-)$-Discrete Dirac Mechanics}
\label{sec:-DDM}
\subsubsection{$(-)$-Discrete Lagrange--Dirac Systems}
Let us first introduce the $(-)$-version of the Dirac differential: Define $\gamma_{Q}^{\rm d-}: T^{*}(Q \times Q) \to T^{*}(Q^{*} \times Q)$ by
\begin{equation*}
  \gamma_{Q}^{\rm d-} \defeq \Omega^{\flat}_{\rm d-} \circ (\kappa_{Q}^{\rm d})^{-1},
\end{equation*}
and, for a given discrete Lagrangian $L_{\rm d}: Q \times Q \to \R$, define the {\rm $(-)$-discrete Dirac differential} $\mathfrak{D}^{-}L_{\rm d}: Q \times Q \to T^{*}(Q^{*} \times Q)$ by
\begin{equation*}
  \mathfrak{D}^{-}L_{\rm d} \defeq \gamma_{Q}^{\rm d-} \circ dL_{\rm d}.
\end{equation*}
In coordinates, we have
\begin{equation*}
  \mathfrak{D}^{-}L_{\rm d}(q^{-}_{k+1}, q_{k+1}) = (-D_{1}L_{\rm d}, q_{k+1}, -q^{-}_{k+1}, -D_{2}L_{\rm d}).
\end{equation*}

\begin{definition}[$(-)$-Discrete Lagrange--Dirac System]
  \label{def:I-DLS}
  Suppose that a discrete Lagrangian $L_{\rm d}: Q \times Q \to \R$ and the constraint distribution $\Delta_{Q} \subset TQ$ are given; and so Eq.~\eqref{eq:Delta_Q^d} gives the discrete constraint distribution $\Delta_{Q}^{\rm d-} \subset Q \times Q$.
  Let $X_{\rm d}^{k}$ be a discrete analogue of a vector field on $T^{*}Q$ as in Eq.~\eqref{eq:X_d}.
  Then, a {\em $(-)$-discrete Lagrange--Dirac system} is a triple $(L_{\rm d}, \Delta_{Q}, X_{\rm d})$ with
  \begin{equation}
    \label{eq:I-DLS}
    \parentheses{ X_{\rm d}^{k}, \mathfrak{D}^{-}L_{\rm d}(q^{-}_{k+1},q_{k+1}) } \in D_{\Delta_{Q}}^{\rm d-}.
  \end{equation}
\end{definition}

\begin{remark}
  The variable $q^{-}_{k+1}$ in Eq.~\eqref{eq:I-DLS} is a discrete analogue of $v$ in Eq.~\eqref{eq:LDS-coord}. See Proposition~\ref{prop:Delta_Q^d}.
\end{remark}

Let us find a coordinate expression for a $(-)$-discrete Lagrange--Dirac system: Eq.~\eqref{eq:I-DLS} gives
\begin{equation*}
  (q_k, q_{k+1}) \in \Delta_{Q}^{\rm d-},
  \qquad
  \mathfrak{D}^{-}L_{\rm d}-\Omega^{\flat}_{\rm d-}(X_{\rm d}^{k}) \in \Delta_{Q^{*} \times Q}^\circ.
\end{equation*}
where
\begin{equation*}
  \Omega^{\flat}_{\rm d-}(X_{\rm d}^{k}) = (p_{k}, q_{k+1}, -q_{k}, -p_{k+1}).
\end{equation*}
Thus, we obtain the following set of equations:
\begin{subequations}
  \label{eq:-DLDEq}
  \begin{equation}
    \begin{array}{c}
      \DS (q_k,q_{k+1}) \in \Delta_{Q}^{\rm d-},
      \qquad
      \DS q_{k} = q^{-}_{k+1},
      \medskip\\
      \DS p_{k} = -D_{1}L_{\rm d}(q^{-}_{k+1}, q_{k+1}),
      \qquad
      \DS p_{k+1} - D_{2}L_{\rm d}(q^{-}_{k+1}, q_{k+1}) \in \Delta_{Q}^\circ(q_{k+1}),
    \end{array}
  \end{equation}
  or more explicitly
  \begin{equation}
    \begin{array}{c}
      \DS \omega_{\rm d-}^{a}(q_k,q_{k+1}) = 0,
      \qquad
      \DS q_{k} = q^{-}_{k+1},
      \medskip\\
      \DS p_{k} = -D_{1}L_{\rm d}(q^{-}_{k+1}, q_{k+1}),
      \qquad
      \DS p_{k+1} - D_{2}L_{\rm d}(q^{-}_{k+1}, q_{k+1}) = \mu_{a} \omega^{a}(q_{k+1}),
    \end{array}
  \end{equation}
\end{subequations}
where $a = 1, 2, \dots, m$.
We shall call them the {\em $(-)$-discrete Lagrange--Dirac equations}; they again recover the nonholonomic integrator of \citet{CoMa2001} (see also \citet{McPe2006}).

If $\Delta_{Q} = TQ$, then $\Delta_{Q}^{\rm d-} = Q\times Q$ and $\Delta_{Q}^{\circ} = 0$; and so the above equations reduce to
\begin{equation}
  q_{k} = q^{-}_{k+1},
  \qquad
  p_{k} = -D_{1}L_{\rm d}(q^{-}_{k+1}, q_{k+1}),
  \qquad
  p_{k+1} = D_{2}L_{\rm d}(q^{-}_{k+1}, q_{k+1}).
\end{equation}
This is a slightly different (but equivalent) expression for Eq.~\eqref{eq:DIELEq}.

\subsubsection{$(-)$-Discrete Nonholonomic Hamiltonian System}
The corresponding discrete nonholonomic Hamiltonian system is defined analogously:
\begin{definition}[$(-)$-Discrete Nonholonomic Hamiltonian System]
  \label{def:I-DHS}
  Suppose that a $(-)$-discrete Hamiltonian~(referred to as the left discrete Hamiltonian in \citep{LaWe2006}) $H_{\rm d-}: Q^{*} \times Q \to \R$ and the constraint distribution $\Delta_{Q} \subset TQ$ are given; and so Eq.~\eqref{eq:Delta_Q^d} gives the discrete constraint distribution $\Delta_{Q}^{\rm d-} \subset Q \times Q$.
  Let $X_{\rm d}^{k}$ be a discrete analogue of a vector field on $T^{*}Q$ as in Eq.~\eqref{eq:X_d}.
  Then, a {\em $(-)$-discrete nonholonomic Hamiltonian system} is a triple $(H_{\rm d-}, \Delta_{Q}, X_{\rm d})$ with
  \begin{equation}
    \label{eq:I-DHS}
    \parentheses{X_{\rm d}^{k}, dH_{\rm d-}(p_{k}, q_{k+1})} \in D_{\Delta_{Q}}^{\rm d-}.
  \end{equation}
\end{definition}

A coordinate expression is obtained in a similar way: We obtain the following set of equations, which we shall call the {\em $(-)$-discrete nonholonomic Hamilton's equations}:
\begin{subequations}
  \label{eq:I-DHS-coord}
  \begin{equation}
    (q_k, q_{k+1}) \in \Delta_{Q}^{\rm d-},
    \qquad
    q_{k} = -D_{1}H_{\rm d-}(p_{k}, q_{k+1}),
    \qquad
    p_{k+1} + D_{2}H_{\rm d-}(p_{k}, q_{k+1}) \in \Delta_{Q}^\circ (q_{k+1}),
  \end{equation}
  or more explicitly
  \begin{equation}
    \omega_{\rm d-}^{a}(q_k,q_{k+1}) = 0,
    \qquad
    q_{k} = -D_{1}H_{\rm d-}(p_{k}, q_{k+1}),
    \qquad
    p_{k+1} + D_{2}H_{\rm d-}(p_{k}, q_{k+1}) = \mu_{a} \omega^{a}(q_{k+1}),
  \end{equation}
\end{subequations}
where $a = 1, 2, \dots, m$.

If $\Delta_{Q} = TQ$, then $\Delta_{Q}^{\rm d-} = Q\times Q$ and $\Delta_{Q}^{\circ} = 0$; and so the above equations reduce to
\begin{equation}
  q_{k} = -D_{1}H_{\rm d-}(p_{k}, q_{k+1}),
  \qquad
  p_{k+1} = -D_{2}H_{\rm d-}(p_{k}, q_{k+1}),
\end{equation}
which are the left discrete Hamilton's equations in Lall and West~\citep{LaWe2006}.

\section{Example of Discrete Lagrange--Dirac System---LC Circuit}
\label{sec:LC_circuit}
\subsection{Formulation}
We apply the above formulation of discrete Dirac mechanics, in particular discrete Lagrange--Dirac systems, to the LC circuit example from Section~\ref{sec:DS_TT_LDS}.

Choose the retraction $\mathcal{R}: TQ \to Q$ (see Section~\ref{ssec:Retractions} for more details) defined by
\begin{equation}
  \label{eq:R-LC}
  \mathcal{R}_{q}(v) \defeq q + v h,
\end{equation}
where $h$ is the time step; hence we have
\begin{equation*}
  \mathcal{R}_{q_{0}}^{-1}(q_{1}) = \frac{q_{1} - q_{0}}{h}.
\end{equation*}
Then, we define the discrete Lagrangian $L_{\rm d}: Q \times Q \to \R$ in terms of the continuous Lagrangian, Eq.~\eqref{eq:L-LC}, as follows:
\begin{align}
  \label{eq:L_d-LC}
  L_{\rm d}(q_{k}, q^{+}_{k})
  &\defeq
  h\,L\parentheses{ q_{k}, \mathcal{R}_{q_{k}}^{-1}(q_{k}^{+}) }
  \nonumber\\
  &=
  h \brackets{
    \frac{\ell}{2} \parentheses{ \frac{q^{+,\ell}_{k} - q^{\ell}_{k}}{h} }^{2}
    - \sum_{i=1}^{3} \frac{(q^{c_{i}}_{k})^{2}}{2 c_{i}}
  }.
\end{align}
This is a discretization that corresponds to the symplectic Euler method~(see, e.g., \cite{MaWe2001}).

We also introduce the discrete constraints $\Delta_{Q}^{\rm d+}$ using Eq.~\eqref{eq:omega_d+-} with the original constraint one-forms $\{\omega^{1}, \omega^{2}\}$ given in Eq.~\eqref{eq:omega-LC}:
\begin{equation*}
  \omega_{\rm d+}^{a}(q_{k}, q_{k+1})
  \defeq
  \ip{ \omega^{a}(q_{k}) }{ \mathcal{R}_{q_{k}}^{-1}(q_{k+1}) }.
\end{equation*}
Simple computations show that
\begin{equation*}
  \begin{array}{l}
    \displaystyle \omega_{\rm d+}^{1}(q_{k}, q_{k+1})
    = \frac{1}{h} \brackets{ -(q^{\ell}_{k+1} - q^{\ell}_{k}) + (q^{c_{2}}_{k+1} - q^{c_{2}}_{k}) },
    \medskip\\
    \displaystyle \omega_{\rm d+}^{2}(q_{k}, q_{k+1})
    = \frac{1}{h} \brackets{ -(q^{c_{1}}_{k+1} - q^{c_{1}}_{k}) + (q^{c_{2}}_{k+1} - q^{c_{2}}_{k}) - (q^{c_{3}}_{k+1} - q^{c_{3}}_{k}) }.
  \end{array}
\end{equation*}
Then, Eq.~\eqref{eq:Delta_Q^d} gives
\begin{align*}
  \Delta_{Q}^{\rm d+}
  &\defeq \setdef{ (q_{k}, q_{k+1}) \in Q \times Q }{ \omega_{\rm d+}^{a}(q_{k}, q_{k+1}) = 0,\, a = 1,2 }
  \nonumber\\
  &= \Bigl\{ (q_{k}, q_{k+1}) \in Q \times Q
  \ |\ {-}q^{\ell}_{k+1} + q^{c_{2}}_{k+1} = -q^{\ell}_{k} + q^{c_{2}}_{k},
  \nonumber\\
  &\hspace{1.75in}
  -q^{c_{1}}_{k+1} + q^{c_{2}}_{k+1} - q^{c_{3}}_{k+1} = -q^{c_{1}}_{k} + q^{c_{2}}_{k} - q^{c_{3}}_{k}
  \Bigr\}.
\end{align*}
Note that the original constraints are holonomic, i.e., the one-forms $\omega^{a}$ are exact, and the above expression for the discrete constraints are the integral form of the original constraints.

Then the $(+)$-discrete Lagrange--Dirac equations~\eqref{eq:+DLDEq} give
\begin{equation}
  \label{eq:I+DLS-LC}
  \begin{array}{c}
    \DS q^{\ell}_{k+1} - q^{\ell}_{k} = q^{c_{2}}_{k+1} - q^{c_{2}}_{k},
    \qquad
    \DS q^{c_{1}}_{k+1} -q^{c_{1}}_{k} = (q^{c_{2}}_{k+1} - q^{c_{2}}_{k}) - (q^{c_{3}}_{k+1} - q^{c_{3}}_{k}),
    \medskip\\
    \DS q^{\ell}_{k+1} = q^{+,\ell}_{k},
    \qquad
    \DS q^{c_{i}}_{k+1} = q^{+,c_{i}}_{k}\quad (i = 1,2,3),
    \medskip\\
    \DS p_{\ell,k+1} = \ell\,\frac{q^{+,\ell}_{k} - q^{\ell}_{k}}{h},
    \qquad
    \DS p_{c_{i},k+1} = 0 \quad (i = 1,2,3),
    \medskip\\
    \DS p_{\ell,k} - \ell\,\frac{q^{+,\ell}_{k} - q^{\ell}_{k}}{h} = -\mu_{1},
    \medskip\\
    \DS p_{c_{i},k} - \frac{h\,q^{+,c_{i}}_{k}}{2 c_{i}} = -\mu_{2}\quad (i = 1,3),
    \qquad
    \DS p_{c_{2},k} - \frac{h\,q^{+,c_{2}}_{k}}{2 c_{2}} = \mu_{1} + \mu_{2},
 \end{array}
\end{equation}
where $\mu_{a}$ are Lagrange multipliers, and we used the fact that $\Delta_{Q}^{\circ} = \Span\{ \omega^{1}, \omega^{2} \}$ with $\omega^{1}$ and $\omega^{2}$ defined in Eq.~\eqref{eq:omega-LC}.

\subsection{Numerical Result}
\label{ssec:NumericalResult}
Assume the initial condition
\begin{equation*}
  q^{\ell}(0) = q^{c_{1}}(0) = q^{c_{2}}(0) = q^{c_{3}}(0) = 0,
  \qquad
  \dot{q}^{\ell}(0) = \dot{q}^{c_{2}}(0) = 10,
  \qquad
  \dot{q}^{c_{1}}(0) = \dot{q}^{c_{3}}(0) = 0.
\end{equation*}
Applying elementary circuit theory to the example, we obtain the exact solution
\begin{equation*}
  q^{\ell}_{\rm ex}(t) = \frac{10}{c} \sin c t,
\end{equation*}
where
\begin{equation*}
  c \defeq \sqrt{\frac{c_{1} + c_{2} + c_{3}}{c_{2}(c_{1} + c_{3}) L}},
\end{equation*}
and thus the period of the solution is $T \defeq 2\pi/c$.
With the choice of the parameters
\begin{equation*}
  \ell = \frac{3}{4},
  \qquad
  c_{1} = 1,
  \qquad
  c_{2} = 2,
  \qquad
  c_{3} = 3,
\end{equation*}
we have $c = 1$, and so the period $T$ becomes $2\pi$.

\begin{figure}[ht!]
  \centering
  \includegraphics[width=.7\linewidth]{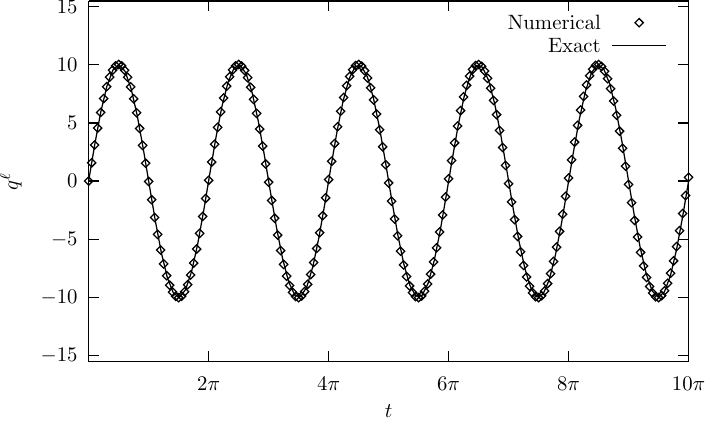}
  \caption{Comparison of exact and numerical solutions (40 points per period) for LC circuit.}
  \label{fig:NumericalResult-LLC}
\end{figure}

Fig.~\ref{fig:NumericalResult-LLC} compares the exact solution with the numerical solution for time step size $h = 2\pi/40 \simeq 0.157$, i.e., 40 time intervals per period.

Table~\ref{tab:Convergence} shows how the error at $t = 5T = 10\pi$ converges as $N$, the number of time intervals per period, increases. The method clearly exhibits second-order convergence behavior, whereas the discretization corresponds to the symplectic Euler method, which is first-order accurate.
\begin{table}[ht!]
  \centering
  \caption{Convergence of numerical method: Number of time intervals per period $N$ vs.~Error at $t = 5T = 10\pi$.}
  \label{tab:Convergence}
  \renewcommand{\arraystretch}{1.25}
  \begin{tabular}{|c||c|c|c|c|}
    \hline
    $N$ & 20 & 40 & 80 & 160
    \\\hline
    $|q^{\ell}_{5N} - q^{\ell}_{\rm ex}(5T)|$ & 1.31915 & 0.324829 & 0.0808631 & 0.0201938
    \\\hline
  \end{tabular}
\end{table}

\begin{remark} One possible explanation for the second-order convergence rate is the following: As one can see from Eq.~\eqref{eq:I+DLS-LC}, the $\{ q^{\ell}_{k} \}$ are the only variables explicitly involved with the time evolution\footnote{This is due to the fact that the original Lagrangian, Eq.~\eqref{eq:L-LC}, is degenerate, i.e., its $f$-dependence is only through $\ell$-component $f^{\ell}$.} and the other variables could be determined from the constraints.
Since $\{ q^{\ell}_{k} \}$ are not present in the potential term in the discrete Lagrangian, Eq.~\eqref{eq:L_d-LC}, only the first term (that corresponds to the inductance energy or ``kinetic energy'' with the electrical-mechanical analogy) is relevant to the time evolution.
However, since the coefficient of this term is constant, the approximation of the ``kinetic energy'' term in the discrete Lagrangian, Eq.~\eqref{eq:L_d-LC}, is the same as that of the midpoint rule, i.e., the approximation given by the discrete Lagrangian of the form
\begin{equation*}
  L_{\rm d}^{\rm MP}(q_{k}, q^{+}_{k})
  =
  h\,L\parentheses{ \frac{q_{k} + q_{k}^{+}}{2}, \frac{q_{k}^{+} - q_{k}}{h} },
\end{equation*}
which yields a second-order accurate method.
\end{remark}

\begin{remark}
  Eliminating $p$ and $\mu$ from Eq.~\eqref{eq:LDS-coord-LC}, we obtain
  \begin{equation*}
    \ell\,\ddot{q}^{\ell} = -\frac{ q^{c_{3}} }{c_{3}} - \frac{ q^{c_{2}} }{c_{2}},
    \qquad
    \dot{q}^{c_{2}} = \dot{q}^{\ell},
    \qquad
    (c_{1} + c_{3})\,\dot{q}^{c_{3}} = c_{3}\,\dot{q}^{c_{2}},
  \end{equation*}
  If we apply the central difference approximation to $\ddot{q}^{\ell}$ and forward difference to all the first-order derivatives in the above equations, we obtain the same numerical method defined by Eq.~\eqref{eq:I+DLS-LC} (after $p_{k}$ and $\mu$ are eliminated).
\end{remark}

In this paper, we do not delve into the issue of accuracy of the numerical methods defined by discrete Lagrange--Dirac systems, instead, we leave it as a topic for future studies.

\section{Variational Structure for Lagrange--Dirac and Nonholonomic Hamiltonian Systems}
\label{sec:VariationalStructure}
In this section we briefly come back to the continuous setting discussed in Section~\ref{sec:DS_TT_LDS} to review variational formulations of Lagrange--Dirac and nonholonomic Hamiltonian systems, again following \citet{YoMa2006b}.
This section is a precursor to the development of the corresponding discrete analogues to follow in the next section.
\subsection{Lagrange--d'Alembert--Pontryagin Principle and Lagrange--Dirac Systems}
\begin{definition}
  \label{def:LdAPP}
  Suppose that a Lagrangian $L: TQ \to \R$ and a constraint distribution $\Delta_{Q} \subset TQ$ are given.
  The {\em Lagrange--d'Alembert--Pontryagin principle} is the augmented variational principle on the \textit{Pontryagin bundle} $TQ\oplus T^*Q$ defined by
  \begin{equation}
    \label{eq:LdAPP}
    \delta \int_{a}^{b} [ L(q,v) + p(\dot{q} - v) ]\,dt = 0,
  \end{equation}
  with the constraint $\dot{q} \in \Delta_{Q}$; we assume that the variation $\delta q$ vanishes at the endpoints, i.e., $\delta q(a) = \delta q(b) = 0$, and also impose $\delta q \in \Delta_{Q}$ {\em after} taking the variations inside the integral sign.
\end{definition}

The Lagrange--Dirac system follows from the Lagrange--d'Alembert--Pontryagin principle:
In a local trivialization, $Q$ is represented by an open set $U$ in a linear space $E$, so the Pontryagin bundle is represented by $(U\times E)\oplus (U\times E^*)\cong U\times E\times E^*$, with local coordinates $(q,v,p)$. If we consider $q$, $v$, and $p$ as independent variables, we have that,
\begin{align*}
  \delta\int_{a}^{b} [L(q,v) + p(\dot{q} - v)]\,dt
  &=\int_{a}^{b} \left[
    \pd{L}{q}\,\delta q
    + \left( \pd{L}{v} - p \right) \delta v
    + (\dot{q} - v)\,\delta p
    +p\,\delta\dot{q}
  \right] dt
  \\
  &=\int_{a}^{b} \left[
    \left( \pd{L}{q} - \dot{p} \right) \delta q
    +\left( \pd{L}{v} - p \right) \delta v
    + (\dot{q} - v)\,\delta p
  \right] dt,
\end{align*}
where we used integration by parts, and the fact that the variation $\delta q$ vanishes at the endpoints.
Taking account of the constraints $\delta q \in \Delta_{Q}$, Eq.~\eqref{eq:LdAPP} gives the Lagrange--Dirac equation~\eqref{eq:LDS-coord}:
\begin{equation}
  \dot{q} \in \Delta_{Q},
  \qquad
  \dot{q} = v,
  \qquad
  p = \pd{L}{v},
  \qquad
  \dot{p} - \pd{L}{q} \in \Delta_{Q}^{\circ}.
\end{equation}

\subsection{Hamilton--d'Alembert Principle in Phase Space and Nonholonomic Hamiltonian Systems}
\begin{definition}
  \label{def:HdAP}
  Suppose that a Hamiltonian $H: T^{*}Q \to \R$ and a constraint distribution $\Delta_{Q} \subset TQ$ are given.
  The {\em Hamilton--d'Alembert principle in phase space} is the variational principle defined by
  \begin{equation}
    \label{eq:HdAP}
    \delta \int_{a}^{b} [ p\,\dot{q} - H(q,p) ]\,dt = 0,
  \end{equation}
  with the constraint $\dot{q} \in \Delta_{Q}$; we assume that the variation $\delta q$ vanishes at the endpoints, i.e., $\delta q(a) = \delta q(b) = 0$, and also impose $\delta q \in \Delta_{Q}$ {\em after} taking the variations inside the integral sign.
\end{definition}

The nonholonomic Hamiltonian system follows from the Hamilton--d'Alembert principle in phase space:
Eq.~\eqref{eq:HdAP} gives
\begin{align*}
  0 = \delta \int_{a}^{b} [ p\,\dot{q} -H(q,p)] dt
  &= \int_{a}^{b} \left(
    \dot{q}\,\delta p
    + p\,\delta\dot{q} - \pd{H}{q}\,\delta q
    - \pd{H}{p}\,\delta p
  \right) dt
  \\
  &= \int_{a}^{b} \left[\left(- \dot p -\frac{\partial H}{\partial q}\right)\delta q+\left(\dot q-\frac{\partial H}{\partial p}\right) \delta p\right]dt,
\end{align*}
which, under the constraints $\delta q \in \Delta_{Q}$, yields
\begin{equation}
  \dot{q} \in \Delta_{Q},
  \qquad
  \dot{q} = \pd{H}{p},
  \qquad
  \dot{p} + \pd{H}{q} \in \Delta_{Q}^{\circ}.
\end{equation}

\section{Discrete Variational Structure for Discrete Lagrange--Dirac and Nonholonomic Hamiltonian Systems}
\label{sec:DiscreteVariationalStructure}
This section develops discrete analogues of the variational structure discussed in the last section.
It is shown that the discrete versions of Lagrange--d'Alembert--Pontryagin principle and Hamilton--d'Alembert principle in phase space yield discrete Lagrange--Dirac and nonholonomic Hamiltonian systems, respectively.
\subsection{Discrete Pontryagin Bundles}
Let us first introduce discrete analogues of the Pontryagin bundle $TQ \oplus T^{*}Q$:
\begin{definition}[$(\pm)$-Discrete Pontryagin Bundles]
  The {\em $(+)$-discrete Pontryagin bundle} is defined by
  \begin{equation*}
    (Q\times Q) \oplus (Q \times Q^{*}) = \braces{ \parentheses{ (q_{k}, q_{k}^{+}), (q_{k}, p_{k+1}) } },
  \end{equation*}
  or, by identifying the first $Q$ of each, we have
  \begin{equation*}
    (Q\times Q) \oplus (Q \times Q^{*}) \cong Q \times Q \times Q^{*} = \braces{ (q_{k}, q_{k}^{+}, p_{k+1}) }.
  \end{equation*}
  Similarly, the {\em $(-)$-discrete Pontryagin bundle} is defined by
  \begin{equation*}
    (Q\times Q) \oplus (Q^{*} \times Q) = \braces{ \parentheses{ (q_{k+1}^{-}, q_{k+1}), (p_{k}, q_{k+1}) } },
  \end{equation*}
  or, by identifying the second $Q$ of each, we have
  \begin{equation*}
    (Q\times Q) \oplus (Q^{*} \times Q) \cong Q \times Q^{*} \times Q = \braces{ (q_{k+1}^{-}, p_{k}, q_{k+1}) }.
  \end{equation*}
\end{definition}

\subsection{Discrete Lagrange--d'Alembert--Pontryagin Principle and Discrete Lagrange--Dirac Systems}
\begin{definition}[$(\pm)$-Discrete Lagrange--d'Alembert--Pontryagin Principle]
  \label{def:DLdAPP}
  Suppose that a discrete Lagrangian $L_{\rm d}: Q \times Q \to \R$ and the constraint distribution $\Delta_{Q} \subset TQ$ are given; and so Eq.~\eqref{eq:Delta_Q^d} gives the discrete constraint distributions $\Delta_{Q}^{\rm d\pm} \subset Q \times Q$.
  Then, the {\em $(\pm)$-discrete Lagrange--d'Alembert--Pontryagin principle} is the discrete augmented variational principle defined by
  \begin{equation}
    \label{eq:+DLdAPP}
    \delta \sum_{k=0}^{N-1} \brackets{L_{\rm d}(q_k, q^{+}_{k}) + p_{k+1}(q_{k+1} - q^{+}_{k})} = 0
  \end{equation}
  or
  \begin{equation}
    \label{eq:-DLdAPP}
    \delta \sum_{k=0}^{N-1} \brackets{L_{\rm d}(q^{-}_{k+1}, q_{k+1}) - p_{k}(q_{k} - q^{-}_{k+1})} = 0,
  \end{equation}
  with the constraint $(q_{k}, q_{k+1}) \in \Delta_{Q}^{\rm d\pm}$ respectively; we assume that the variations $\delta q_{k}$ vanish at the endpoints, i.e., $\delta q_{0} = \delta q_{N} = 0$, and also impose $\delta q_{k} \in \Delta_{Q}(q_{k})$ {\em after} taking the variations inside the summation.
\end{definition}

\begin{proposition}
  The $(\pm)$-discrete Lagrange--d'Alembert--Pontryagin principles yield the $(\pm)$-discrete Lagrange--Dirac equations~\eqref{eq:+DLDEq} and \eqref{eq:-DLDEq}, respectively.
\end{proposition}

\begin{proof}
  First taking the variations in Eqs.~\eqref{eq:+DLdAPP} and \eqref{eq:-DLdAPP}, we have
  \begin{align*}
    0 &= \delta\sum_{k=0}^{N-1} \brackets{L_{\rm d}(q_{k}, q^{+}_{k}) + p_{k+1}(q_{k+1} - q^{+}_{k})}
    \\
    &= \sum_{k=1}^{N-1} \brackets{ D_{1}L_{\rm d}(q_{k}, q^{+}_{k})+p_{k} }\delta q_{k}
    + \sum_{k=0}^{N-1} \braces{
      \brackets{
        D_{2}L_{\rm d}(q_{k},q^{+}_{k})-p_{k+1}
      }\delta q^{+}_{k}
      + \parentheses{ q_{k+1} - q^{+}_{k} }\delta p_{k+1}
    },
\intertext{and}
    0 &= \delta\sum_{k=0}^{N-1} \brackets{L_{\rm d}(q^{-}_{k+1}, q_{k+1}) - p_{k}(q_{k} - q^{-}_{k+1})}
    \\
    &= \sum_{k=0}^{N-2} \brackets{ D_{2}L_{\rm d}(q^{-}_{k+1}, q_{k+1}) - p_{k+1} }\delta q_{k+1}\\
    &\qquad\qquad\qquad +\sum_{k=0}^{N-1} \braces{
      \brackets{ D_{1}L_{\rm d}(q^{-}_{k+1}, q_{k+1}) + p_{k} }\delta q^{-}_{k+1}
      + \parentheses{ q_{k+1}^{-} - q_{k} }\delta p_{k}
    },
  \end{align*}
  where we used $\delta q_{0} = 0$ and $\delta q_{N} = 0$.
  Taking account of the corresponding constraints on the variations in each of the above equations, we obtain Eqs.~\eqref{eq:+DLDEq} and \eqref{eq:-DLDEq}, respectively.
\end{proof}

\subsection{Discrete Hamilton--d'Alembert Principle in Phase Space and Discrete Nonholonomic Hamiltonian Systems}
\begin{definition}[$(\pm)$-Discrete Hamilton--d'Alembert Principle in Phase Space]
  \label{def:DHdAP}
  Suppose that a $(\pm)$-discrete Hamiltonian $H_{\rm d+}: Q \times Q^{*} \to \R$ or $H_{\rm d-}: Q^{*} \times Q \to \R$ and the constraint distribution $\Delta_{Q} \subset TQ$ are given; and so Eq.~\eqref{eq:Delta_Q^d} gives the discrete constraint distributions $\Delta_{Q}^{\rm d\pm} \subset Q \times Q$.
  Then, the {\em $(\pm)$-discrete Hamilton--d'Alembert principle in phase space} is the discrete variational principle defined by
  \begin{equation}
    \label{eq:D+HdAP}
    \delta \sum_{k=0}^{N-1} [p_{k+1} q_{k+1} - H_{\rm d+}(q_k, p_{k+1})] = 0
  \end{equation}
  or
  \begin{equation}
    \label{eq:D-HdAP}
    \delta \sum_{k=0}^{N-1} [-p_{k} q_{k} - H_{\rm d-}(p_{k}, q_{k+1})] = 0,
  \end{equation}
  with the constraint $(q_{k}, q_{k+1}) \in \Delta_{Q}^{\rm d\pm}$ respectively; we assume that the variations $\delta q_{k}$ vanish at the endpoints, i.e., $\delta q_{0} = \delta q_{N} = 0$, and also impose $\delta q_{k} \in \Delta_{Q}(q_{k})$ {\em after} taking the variations inside the summation.
\end{definition}

\begin{proposition}
  The $(\pm)$-discrete Hamilton--d'Alembert principles yield the $(\pm)$-discrete nonholonomic Hamilton's equations~\eqref{eq:I+DHS-coord} and \eqref{eq:I-DHS-coord}, respectively.
\end{proposition}

\begin{proof}
  First taking the variations in Eqs.~\eqref{eq:D+HdAP} and \eqref{eq:D-HdAP}, we have
  \begin{align*}
    0 &= \delta \sum_{k=0}^{N-1} [ p_{k+1} q_{k+1}-H_{\rm d+}(q_k,p_{k+1})]
    \\
    &= \sum_{k=0}^{N-1} \brackets{ q_{k+1} - D_{2}H_{\rm d+}(q_k, p_{k+1}) } \delta p_{k+1}
    + \sum_{k=1}^{N-1} \brackets{ p_{k} - D_{1}H_{\rm d+}(q_k,p_{k+1}) } \delta q_{k}
\intertext{and}
    0 &= \delta \sum_{k=0}^{N-1} [ -p_{k} q_{k} - H_{\rm d-}(p_k, q_{k+1})]
    \\
    &= -\sum_{k=0}^{N-1} \brackets{ q_{k} + D_{1}H_{\rm d-}(p_{k}, q_{k+1}) } \delta p_{k}
    - \sum_{k=0}^{N-2} \brackets{ p_{k+1} + D_{2}H_{\rm d-}(p_{k}, q_{k+1}) } \delta q_{k+1},
  \end{align*}
  where we used $\delta q_{0} = 0$ and $\delta q_{N} = 0$.
  Taking account of the constraints on the variations $\delta q_{k}$ in each of the above equations, we obtain Eqs.~\eqref{eq:I+DHS-coord} and \eqref{eq:I-DHS-coord}, respectively.
\end{proof}

\section{Extension to Computations on Manifolds}
\label{sec:ExtensionToManifolds}
This section presents a means to apply the preceding theory to computations for the case when $Q$ is a manifold.
We do not attempt a full extension of the theory to manifolds since discrete Hamiltonian mechanics \citep{LaWe2006} is not intrinsic: Recall that the $(+)$-discrete Hamiltonian $H_{\rm d+}$ is a Type~2 generating function (see \citet{LaWe2006} and also Section~\ref{ssec:DiscreteMechAndGenFunctions}), which is based on the idea of generating the pair $(p_{0}, q_{1})$ with the pair $(q_{0}, p_{1})$ fixed.
However, this does not make intrinsic sense, since fixing $p_{1}$ in $T^{*}Q$ requires that its corresponding base point $q_{1}$ is fixed as well.

Instead, we make use of the idea of retractions and introduce the notion of retraction compatible coordinate charts to provide a means of applying the results in the linear theory to computations on manifolds, in a semi-globally compatible fashion. Retraction compatible coordinate charts provide a generalization of the canonical coordinates of the first kind on Lie groups (see, e.g., \citet[Section~2.10]{Va1984} and also Example~\ref{ex:RetractionOnLieGroups} below) to more general configuration manifolds. By semi-global, we mean that the discrete flow is well-defined on a neighborhood of the diagonal of $Q\times Q$, which corresponds to a restriction on the size of the time step.

In particular, on a retraction compatible coordinate chart, the discrete flow is described, in local coordinates, by the vector space expressions. This has non-trivial implications for geometric numerical integration, since na\"\i vely applying a linear space numerical integrator on different charts may lead to poor global properties, as discussed in \citep{BeChFa2001}. By restricting ourselves to retraction compatible coordinate charts, we ensure that the local conservation properties of the geometric numerical integrators we introduce in this paper persist globally as well.

\subsection{Retractions}
\label{ssec:Retractions}
Let us first recall the definition of a retraction:
\begin{definition}[{\citet[][Definition~4.1.1 on p.~55]{AbMaSe2008}}]
  A retraction on a manifold $Q$ is a smooth mapping $\mathcal{R}: TQ \to Q$ with the following properties:
  Let $\mathcal{R}_{q}: T_{q}Q \to Q$ be the restriction of $\mathcal{R}$ to $T_{q}Q$ for an arbitrary $q \in Q$; then,
  \begin{enumerate}[(i)]
  \item $\mathcal{R}_{q}(0_{q}) = q$, where $0_{q}$ denotes the zero element of $T_{q}Q$;
    \label{item:Retraction-i}
  \item with the identification $T_{0_{q}}T_{q}Q \simeq T_{q}Q$, $\mathcal{R}_{q}$ satisfies
    \begin{equation}
      \label{eq:TR}
      T_{0_{q}}\mathcal{R}_{q} = \id_{T_{q}Q},
    \end{equation}
    where $T_{0_{q}}\mathcal{R}_{q}$ is the tangent map of $\mathcal{R}_{q}$ at $0_{q} \in T_{q}Q$.
    \label{item:Retraction-ii}
  \end{enumerate}
\end{definition}

\begin{remark}
  \label{rem:mathcalR_inverse}
  Eq.~\eqref{eq:TR} implies that the map $\mathcal{R}_{q}: T_{q}Q \to Q$ is invertible in some neighborhood of $0_{q}$ in $T_{q}Q$.
\end{remark}

It is convenient to introduce $\tilde{\mathcal{R}}: TQ \to Q \times Q$ defined by
\begin{equation}
  \label{eq:tildeR}
  \tilde{\mathcal{R}}(v_{q}) \defeq (q, \mathcal{R}_{q}(v_{q})).
\end{equation}
It is easy to see from the above expression and the above remark that $\tilde{\mathcal{R}}: TQ \to Q \times Q$ is also invertible in some neighborhood of $0_{q} \in TQ$ for any $q \in Q$.
 
Let us introduce a special class of coordinate charts that are convenient to work with:
\begin{definition}[Retraction compatible coordinate charts and atlas]
  Let $Q$ be an $n$-dimensional manifold equipped with a retraction $\mathcal{R}: TQ \to Q$.
  A coordinate chart $(U, \varphi)$ with $U$ an open subset in $Q$ and $\varphi: U \to \R^{n}$ is said to be {\em retraction compatible at} $q \in U$ if
  \begin{enumerate}[(i)]
  \item $\varphi$ is centered at $q$, i.e., $\varphi(q) = 0$;
  \item the compatibility condition
    \begin{equation}
      \label{eq:R-compatibility}
      \mathcal{R}(v_{q}) = \varphi^{-1} \circ T_{q}\varphi(v_{q})
    \end{equation}
    holds, where we identify $T_{0}\R^{n}$ with $\R^{n}$ as follows: Let $\varphi = (x^{1}, \dots, x^{n})$ with $x^{i}: U \to \R$ for $i = 1, \dots, n$. Then
    \begin{equation}
      \label{eq:f_varphi}
      v^{i}\pd{}{x^{i}} \mapsto (v^{1}, \dots, v^{n}),
    \end{equation}
    where $\tpd{}{x^{i}}$ is the unit vector in the $x^{i}$-direction in $T_{0}\R^{n}$.
  \end{enumerate}
  An atlas for the manifold $Q$ is {\em retraction compatible} if it consists of retraction compatible coordinate charts.
\end{definition}

\begin{remark}
  In Eq.~\eqref{eq:R-compatibility}, we assumed that $T_{q}\varphi(v_{q}) \in \varphi(U) \subset \R^{n}$ and so strictly speaking $\mathcal{R}_{q}$ is defined on $(T_{q}\varphi)^{-1}(\varphi(U)) \subset T_{q}Q$.
  However, it is always possible to define a coordinate chart such that $\varphi(U) = \R^{n}$ by ``stretching out'' the open set $\varphi(U)$ to $\R^{n}$ so that Eq.~\eqref{eq:R-compatibility} is defined for any $v_{q} \in T_{q}Q$.
\end{remark}

\begin{example}[Retraction and canonical coordinates of the first kind on a Lie group]
  \label{ex:RetractionOnLieGroups}
  Let $G$ be a (finite-dimensional) Lie group and $\mathfrak{g}$ be its Lie algebra.
  The exponential map $\exp: \mathfrak{g} \to G$ (see, e.g., \citet[Section~9.1]{MaRa1999} and \citet[Section~2.10]{Va1984}) is a diffeomorphism on an open neighborhood $\mathfrak{u}$ of the origin of $\mathfrak{g}$.
  Let $U$ be the neighborhood of the identity $e$ in $G$ defined by $U \defeq \exp(\mathfrak{u}) \subset G$, and restrict the domain of the exponential map to redefine $\exp: \mathfrak{u} \to U$ for notational simplicity.
  Then, it is a diffeomorphism and so we have the inverse $\exp^{-1}: U \to \mathfrak{u}$.

  Now let us define $\mathcal{R}_{g}: T_{g}G \to G$ for any $g \in G$ by\footnote{Strictly speaking, $\mathcal{R}_{g}$ is defined only on $T_{e}L_{g}(\mathfrak{u}) \subset T_{g}G$.}
  \begin{equation*}
    \mathcal{R}_{g} \defeq L_{g} \circ \exp \circ\, T_{g}L_{g^{-1}},
  \end{equation*}
  where $L_{g}: G \to G$ is the left translation by $g$.
  This indeed gives a retraction: Since $\exp(0) = e$, we have $\mathcal{R}_{g}(0_{g}) = g$; we also have, with the identification $T_{0_{g}}T_{g}G \simeq T_{g}G$, 
  \begin{align*}
    T_{0_{g}}\mathcal{R}_{g}
    &= T_{e}L_{g} \circ T_{0}\exp \circ T_{0_{g}}T_{g}L_{g^{-1}}
    \\
    &= T_{e}L_{g} \circ T_{g}L_{g^{-1}}
    \\
    &= \id_{T_{g}G},
  \end{align*}
  where we used the fact that $T_{0}\exp: T\mathfrak{u} \simeq \mathfrak{g} \to \mathfrak{g}$ is the identity (see \cite[Eq.~(2.10.17) on p.~88]{Va1984}), and also that $T_{0_{g}}T_{g}L_{g^{-1}} = T_{g}L_{g^{-1}}$ with the above identification\footnote{The derivative of a linear map at the origin is the linear map itself.}.

  The exponential map also induces the {\em canonical coordinates of the first kind} on the Lie group $G$ as follows~(see, e.g., \citet[Section~2.10]{Va1984} and \citet{MaPeSh1999}): For any $g \in G$, let $U_{g} \defeq L_{g}(U)$ and define a chart $\varphi_{g}: U_{g} \to \mathfrak{g}$ by
  \begin{equation*}
    \varphi_{g} \defeq \exp^{-1} \circ L_{g^{-1}}.
  \end{equation*}
  Then, the chart $\varphi_{g}$ is retraction compatible: We have
  \begin{equation*}
    \varphi_{g}(g) = \exp^{-1} \circ L_{g^{-1}}(g) = \exp^{-1}(e) = 0,
  \end{equation*}
  and also, with the identification $T\mathfrak{u} \simeq \mathfrak{g}$, 
  \begin{align*}
    \varphi_{g}^{-1} \circ T_{g}\varphi_{g}
    &= L_{g} \circ \exp \circ\, T_{e}\exp^{-1} \circ\, T_{g}L_{g^{-1}}
    \\
    &= L_{g} \circ \exp \circ\, T_{g}L_{g^{-1}}
    \\
    &= \mathcal{R}_{g},
  \end{align*}
  where we used the fact that $T_{e}\exp^{-1} = \id_{\mathfrak{g}}$, which follows from $T_{0}\exp = \id_{\mathfrak{g}}$ mentioned above.
\end{example}

Calculations involving a retraction are particularly simple with a retraction compatible chart:
\begin{proposition}
  \label{prop:pq-pairing}
  Let $(U, \varphi)$ be a retraction compatible chart at a point $q \in U$.
  Take an arbitrary point $r$ in $U$ and let $(r^{1}, \dots, r^{n}) \defeq \varphi(r) \in \R^{n}$.
  Then
  \begin{equation}
    \label{eq:Rinv}
    \mathcal{R}_{q}^{-1}(r) = r^{i} \left.\!\pd{}{x^{i}}\right|_{q}
  \end{equation}
  where
  \begin{equation*}
    \left.\!\pd{}{x^{i}}\right|_{q} \defeq T_{0}\varphi^{-1}\parentheses{ \pd{}{x^{i}} } \in T_{q}Q.
  \end{equation*}
  Furthermore, let $dx^{i}|_{q} \in T_{q}^{*}Q$ be the dual basis to $\tpd{}{x^{i}}|_{q} \in T_{q}Q$, i.e., $dx^{i}|_{q}(\tpd{}{x^{j}}|_{q}) = \delta^{i}_{j}$.
  Then, for any $p_{q} = p_{i}\,dx^{i}|_{q} \in T_{q}^{*}Q$, we have
  \begin{equation}
    \label{eq:p-q_pairing}
    \ip{ p_{q} }{ \mathcal{R}_{q}^{-1}(r) } = \ip{ p_{q} }{ \tilde{\mathcal{R}}^{-1}(q, r) } = p_{i} r^{i},
  \end{equation}
  where $\ip{\,\cdot\,}{\,\cdot\,}$ is the natural pairing between elements in $T^{*}Q$ and $TQ$.
\end{proposition}

\begin{proof}
  Follows from straightforward calculations.
\end{proof}

\subsection{Discrete Lagrange--d'Alembert--Pontryagin Principles with Retraction}
\label{ssec:DHPwithRetraction}
Let us use a retraction to reformulate the {\em $(\pm)$-discrete Lagrange--d'Alembert--Pontryagin principle}, Eq.~\eqref{eq:+DLdAPP}, as follows:
\begin{definition}[$(\pm)$-Discrete Lagrange--d'Alembert--Pontryagin Principle with Retraction]
  \label{def:DLdAPP-Retraction}
  Suppose that a discrete Lagrangian $L_{\rm d}: Q \times Q \to \R$ and the constraint distribution $\Delta_{Q} \subset TQ$ are given; and so Eq.~\eqref{eq:Delta_Q^d} gives the discrete constraint distributions $\Delta_{Q}^{\rm d\pm} \subset Q \times Q$.
  Then, the {\em $(\pm)$-discrete Lagrange--d'Alembert--Pontryagin principle} is the discrete augmented variational principle defined by
  \begin{equation}
    \label{eq:+DLdAPP-Retraction}
    \delta S_{\rm d+}^{N} =
    \delta \sum_{k=0}^{N-1} \brackets{
      L_{\rm d}(q_{k}, q^{+}_{k}) + \ip{ p_{k+1} }{ \mathcal{R}_{q_{k+1}}^{-1}(q_{k+1}) - \mathcal{R}_{q_{k+1}}^{-1}(q^{+}_{k}) }
    }
  \end{equation}
  or
  \begin{equation}
    \label{eq:-DLdAPP-Retraction}
    \delta S_{\rm d-}^{N} =
    \delta \sum_{k=0}^{N-1} \brackets{
      L_{\rm d}(q_{k+1}^{-}, q_{k+1}) - \ip{ p_{k} }{ \mathcal{R}_{q_{k}}^{-1}(q_{k}) - \mathcal{R}_{q_{k}}^{-1}(q^{-}_{k+1}) }
    },
  \end{equation}
  with the constraint $(q_{k}, q_{k+1}) \in \Delta_{Q}^{\rm d\pm}$ respectively; we assume that the variations $\delta q_{k}$ vanish at the endpoints, i.e., $\delta q_{0} = \delta q_{N} = 0$, and also impose $\delta q_{k} \in \Delta_{Q}(q_{k})$ {\em after} taking the variations inside the summation.
\end{definition}

With a retraction compatible coordinate chart, Lemma~\ref{prop:pq-pairing} implies that Eqs.~\eqref{eq:+DLdAPP-Retraction} and \eqref{eq:-DLdAPP-Retraction} become
\begin{equation*}
  S_{\rm d+}^{N} = \sum_{k=0}^{N-1} \brackets{
    L_{\rm d}(q_{k}, q^{+}_{k}) + p_{k+1} \cdot (q_{k+1} - q^{+}_{k})
  }
\end{equation*}
and
\begin{equation*}
  S_{\rm d-}^{N} = \sum_{k=0}^{N-1} \brackets{
    L_{\rm d}(q_{k+1}^{-}, q_{k+1}) - p_{k} \cdot (q_{k} - q^{-}_{k+1})
  },
\end{equation*}
where we slightly abused the notation, i.e., $q_{k+1}$, $q^{+}_{k}$, $q^{-}_{k+1}$ are interpreted as both points in $Q$ as well as their coordinate representations.
Therefore the $(\pm)$-discrete Lagrange--d'Alembert--Pontryagin principle in Definition~\ref{def:DLdAPP-Retraction}, written in terms of retraction compatible charts, reduce to those in the linear theory, i.e., Definition~\ref{def:DLdAPP}.
\begin{remark}
  Note that $\mathcal{R}_{q_{k}}^{-1}(q_{k}) = 0$ by definition, and so the terms of the form
  \begin{equation*}
    \ip{ p_{k} }{ \mathcal{R}_{q_{k}}^{-1}(q_{k}) } = p_{k} \cdot q_{k}
  \end{equation*}
  vanish.
\end{remark}

The above discussion implies that the discrete Lagrange--Dirac equations~\eqref{eq:+DLDEq} and \eqref{eq:-DLDEq} in the linear theory are coordinate representations (using a retraction compatible chart) of the systems defined by the discrete Lagrange--d'Alembert--Pontryagin principles in Definition~\ref{def:DLdAPP-Retraction}.
Therefore, if we have a retraction compatible atlas on $Q$, then we may use the coordinate expression, Eq.~\eqref{eq:+DLDEq}, in the linear theory to perform a computation in a single chart, and, if necessary, transform the system to another chart in the atlas, which again has the same form as Eq.~\eqref{eq:+DLDEq}, to continue the computation.

\section{Conclusion}
In this paper, we developed the theoretical foundations of discrete Dirac mechanics from two different perspectives: One through discrete analogues of Tulczyjew's triple and induced Dirac structures, and the other from the variational point of view.

We exploited the discrete Tulczyjew triples to define discrete analogues of the Dirac differential, which is a key in defining the discrete Lagrange--Dirac systems, particularly those with degenerate Lagrangians; and we employed the discrete induced Dirac structures to incorporate discrete constraints.
We also introduced extended discrete variational principles, i.e., the discrete Lagrange--d'Alembert--Pontryagin and Hamilton--d'Alembert principles that give variational formulations of discrete Lagrange--Dirac and nonholonomic Hamiltonian systems.

An LC circuit is taken as an example of a system with a degenerate Lagrangian and constraints, and is modeled as a discrete Lagrange--Dirac system.
We performed numerical computations with the resulting scheme and obtained numerical solutions that converge to the exact solution obtained by elementary circuit theory.

Several interesting topics for future work are suggested by the theoretical developments introduced in this paper:
\begin{itemize}
\item {\em Application to inter-connected systems.} Port-Hamiltonian systems~\citep{Sc2006} provide a natural description of modular and interconnected systems, but this does not naturally lead to geometric structure-preserving discretizations of interconnected systems. It is therefore desirable to develop a unified port-Lagrangian framework for modeling and simulating interconnected systems based on extensions of Lagrange--Dirac mechanics and variational discrete Dirac mechanics.
  \smallskip
\item {\em Hamilton--Jacobi theory for Lagrange--Dirac systems}~(\citet{DiracHJ}).
  Since Dirac structures are related to Lagrangian submanifolds, which in turn describe the geometry of the Hamilton--Jacobi equation, it is natural to explore the Dirac description of Hamilton--Jacobi theory.
  The resulting theory is expected to give insights into discrete Dirac mechanics as the classical Hamilton--Jacobi theory does to discrete mechanics~\citep[Sections~1.8 and 4.8]{MaWe2001}; it is also natural to expect it to specialize to nonholonomic Hamilton--Jacobi theory~\citep{IgLeMa2008, LeMaMa2010, OhBl2009, CaGrMaMaMuRo2010, OhFeBlZe2011}.
  \smallskip
\item {\em Discrete reduction theory for discrete Dirac mechanics with symmetry.}
  The Dirac formulation of reduction (see \citet{YoMa2007,YoMa2009}) provides a means of unifying symplectic, Poisson, nonholonomic, Lagrangian, and Hamiltonian reduction theory, as well as addressing the issue of reduction by stages.
  The discrete analogue of Dirac reduction will proceed by considering the issue of quotient discrete Dirac structures, and constructing a category containing discrete Dirac structures, that is closed under quotients.
  \smallskip
\item {\em Discrete multi-Dirac mechanics for Hamiltonian partial differential equations.}
  Dirac generalizations of multisymplectic field theory (see \citet{VaYoMa2010}), and their corresponding discretizations will provide important insights into the construction of geometric numerical methods for degenerate field theories, such as the Einstein equations of general relativity.
  \smallskip
\item {\em Variational error analysis of discrete Lagrange--Dirac systems.}
 It is natural, and desirable, to extend the variational error analysis techniques developed by \citet{MaWe2001} for discrete Lagrangian mechanics to the case of discrete Lagrange--Dirac systems. In particular, this may provide insight into the rather unexpected convergence behavior observed in Section~\ref{ssec:NumericalResult}.
\end{itemize}

\section*{Acknowledgements}
We gratefully acknowledge helpful comments and suggestions of the referees, Henry Jacobs, Jerrold Marsden, Joris Vankerschaver, Hiroaki Yoshimura, and also the reviewer of our earlier work~\citep{LeOh2010}.
This material is based upon work supported by the National Science Foundation under the applied mathematics grant DMS-0726263 and the Faculty Early Career Development (CAREER) award DMS-1010687. 

\bibliography{DiscreteDirac}

\begin{thebibliography}{54}
\providecommand{\natexlab}[1]{#1}
\providecommand{\url}[1]{\texttt{#1}}
\expandafter\ifx\csname urlstyle\endcsname\relax
  \providecommand{\doi}[1]{doi: #1}\else
  \providecommand{\doi}{doi: \begingroup \urlstyle{rm}\Url}\fi

\bibitem[Abraham and Marsden(1978)]{AbMa1978}
R.~Abraham and J.~E. Marsden.
\newblock \emph{Foundations of Mechanics}.
\newblock Addison--Wesley, 2nd edition, 1978.

\bibitem[Absil et~al.(2008)Absil, Mahony, and Sepulchre]{AbMaSe2008}
P.-A. Absil, R.~Mahony, and R.~Sepulchre.
\newblock \emph{Optimization Algorithms on Matrix Manifolds}.
\newblock Princeton University Press, 2008.

\bibitem[Arnold(1989)]{Ar1989}
V.~I. Arnold.
\newblock \emph{Mathematical Methods of Classical Mechanics}.
\newblock Springer, 1989.

\bibitem[Bates and Sniatycki(1993)]{BaSn1993}
L.~Bates and J.~Sniatycki.
\newblock Nonholonomic reduction.
\newblock \emph{Reports on Mathematical Physics}, 32\penalty0 (1):\penalty0
  99--115, 1993.

\bibitem[Benettin et~al.(2001)Benettin, Cherubini, and Fass{\`o}]{BeChFa2001}
G.~Benettin, {A. M.} Cherubini, and F.~Fass{\`o}.
\newblock A changing-chart symplectic algorithm for rigid bodies and other
  {H}amiltonian systems on manifolds.
\newblock \emph{SIAM J. Sci. Comput.}, 23\penalty0 (4):\penalty0 1189--1203,
  2001.
\newblock ISSN 1064-8275.

\bibitem[Bloch(2003)]{Bl2003}
A.~M. Bloch.
\newblock \emph{Nonholonomic Mechanics and Control}.
\newblock Springer, 2003.

\bibitem[Bloch and Crouch(1997)]{BlCr1997}
A.~M. Bloch and P.~E. Crouch.
\newblock Representations of {D}irac structures on vector spaces and nonlinear
  {L-C} circuits.
\newblock In \emph{Differential Geometry and Control Theory}, pages 103--117.
  American Mathematical Society, 1997.

\bibitem[Bou-Rabee and Marsden(2008)]{BoMa2008}
N.~Bou-Rabee and {J. E.} Marsden.
\newblock {H}amilton--{P}ontryagin integrators on {L}ie groups part {I}:
  Introduction and structure-preserving properties.
\newblock \emph{Found. Comput. Math.}, 2008.

\bibitem[Cari\~nena et~al.(2010)Cari\~nena, Gracia, Marmo, Mart\'inez, Mun\~oz
  Lecanda, and Rom\'an-Roy]{CaGrMaMaMuRo2010}
J.~F. Cari\~nena, X.~Gracia, G.~Marmo, E.~Mart\'inez, M.~C. Mun\~oz Lecanda,
  and N.~Rom\'an-Roy.
\newblock Geometric {H}amilton--{J}acobi theory for nonholonomic dynamical
  systems.
\newblock \emph{International Journal of Geometric Methods in Modern Physics},
  7\penalty0 (3):\penalty0 431--454, 2010.

\bibitem[Cervera et~al.(2003)Cervera, van~der Schaft, and
  Ba{\~n}os]{CeScBa2003}
J.~Cervera, A.~J. van~der Schaft, and A.~Ba{\~n}os.
\newblock On composition of {D}irac structures and its implications for control
  by interconnection.
\newblock In \emph{Nonlinear and adaptive control}, volume 281 of \emph{Lecture
  Notes in Control and Inform. Sci.}, pages 55--63. Springer, Berlin, 2003.

\bibitem[Cort\'es and Mart\'inez(2001)]{CoMa2001}
J.~Cort\'es and S.~Mart\'inez.
\newblock Non-holonomic integrators.
\newblock \emph{Nonlinearity}, 14\penalty0 (5):\penalty0 1365--1392, 2001.

\bibitem[Courant(1990{\natexlab{a}})]{Co1990a}
T.~Courant.
\newblock Dirac manifolds.
\newblock \emph{Transactions of the American Mathematical Society},
  319\penalty0 (2):\penalty0 631--661, 1990{\natexlab{a}}.

\bibitem[Courant(1990{\natexlab{b}})]{Co1990b}
T.~Courant.
\newblock Tangent {D}irac structures.
\newblock \emph{Journal of Physics A: Mathematical and General}, 23\penalty0
  (22):\penalty0 5153--5168, 1990{\natexlab{b}}.

\bibitem[Dalsmo and van~der Schaft(1998)]{DaSc1998}
M.~Dalsmo and A.~J. van~der Schaft.
\newblock On representations and integrability of mathematical structures in
  energy-conserving physical systems.
\newblock \emph{SIAM Journal on Control and Optimization}, 37\penalty0
  (1):\penalty0 54--91, 1998.

\bibitem[de~Le\'on et~al.(2010)de~Le\'on, Marrero, and Mart\'in~de
  Diego]{LeMaMa2010}
M.~de~Le\'on, J.~C. Marrero, and D.~Mart\'in~de Diego.
\newblock Linear almost {P}oisson structures and {H}amilton--{J}acobi equation.
  {A}pplications to nonholonomic mechanics.
\newblock \emph{Journal of Geometric Mechanics}, 2\penalty0 (2):\penalty0
  159--198, 2010.

\bibitem[Dirac(1950)]{Di1950}
P.~A.~M. Dirac.
\newblock Generalized {H}amiltonian dynamics.
\newblock \emph{Canad. J. Math.}, 2:\penalty0 129--148, 1950.

\bibitem[Dirac(1958)]{Di1958a}
P.~A.~M. Dirac.
\newblock Generalized {H}amiltonian dynamics.
\newblock \emph{Proceedings of the Royal Society of London. Series A,
  Mathematical and Physical Sciences}, 246\penalty0 (1246):\penalty0 326--332,
  1958.

\bibitem[Dirac(1964)]{Di1964}
P.~A.~M. Dirac.
\newblock \emph{Lectures on quantum mechanics}.
\newblock Belfer Graduate School of Science, Yeshiva University, New York,
  1964.

\bibitem[Goldstein et~al.(2001)Goldstein, Poole, and Safko]{GoPoSa2001}
H.~Goldstein, C.~P. Poole, and J.~L. Safko.
\newblock \emph{Classical Mechanics}.
\newblock Addison Wesley, 3rd edition, 2001.

\bibitem[Gotay and Nester(1979{\natexlab{a}})]{GoNe1979a}
M.~J. Gotay and J.~M. Nester.
\newblock Presymplectic {H}amilton and {L}agrange systems, gauge
  transformations and the {D}irac theory of constraints.
\newblock In \emph{Group Theoretical Methods in Physics}, volume~94, pages
  272--279. Springer, 1979{\natexlab{a}}.

\bibitem[Gotay and Nester(1979{\natexlab{b}})]{GoNe1979b}
M.~J. Gotay and J.~M. Nester.
\newblock Presymplectic {L}agrangian systems. {I}: the constraint algorithm and
  the equivalence theorm.
\newblock \emph{Annales de l'institut Henri Poincar\'e (A)}, 30\penalty0
  (2):\penalty0 129--142, 1979{\natexlab{b}}.

\bibitem[Gotay and Nester(1980)]{GoNe1980}
M.~J. Gotay and J.~M. Nester.
\newblock Presymplectic {L}agrangian systems. {II}: the second-order equation
  problem.
\newblock \emph{Annales de l'institut Henri Poincar\'e (A)}, 32\penalty0
  (1):\penalty0 1--13, 1980.

\bibitem[Hairer et~al.(2006)Hairer, Lubich, and Wanner]{HaLuWa2006}
E.~Hairer, C.~Lubich, and G.~Wanner.
\newblock \emph{Geometric Numerical Integration: Structure-Preserving
  Algorithms for Ordinary Differential Equations}.
\newblock Springer, Berlin, Heidelberg, 2006.

\bibitem[Iglesias et~al.(2008)Iglesias, Marrero, Mart\'in~de Diego, and
  Mart\'inez]{IgMaMaMa2008}
D.~Iglesias, J.~C. Marrero, D.~Mart\'in~de Diego, and E.~Mart\'inez.
\newblock Discrete nonholonomic lagrangian systems on lie groupoids.
\newblock \emph{Journal of Nonlinear Science}, 18\penalty0 (3):\penalty0
  221--276, 2008.

\bibitem[Iglesias-Ponte et~al.(2008)Iglesias-Ponte, de~Le\'on, and Mart\'in~de
  Diego]{IgLeMa2008}
D.~Iglesias-Ponte, M.~de~Le\'on, and D.~Mart\'in~de Diego.
\newblock Towards a {H}amilton--{J}acobi theory for nonholonomic mechanical
  systems.
\newblock \emph{Journal of Physics A: Mathematical and Theoretical},
  41\penalty0 (1), 2008.

\bibitem[Kharevych et~al.(2006)Kharevych, Yang, Tong, Kanso, Marsden,
  Schr\"oder, and Desbrun]{Kh2006}
L.~Kharevych, W.~Yang, Y.~Tong, E.~Kanso, J.~E. Marsden, P.~Schr\"oder, and
  M.~Desbrun.
\newblock Geometric, variational integrators for computer animation.
\newblock In \emph{ACM/EG Symposium on Computer Animation}, pages 43--51, 2006.

\bibitem[Koon and Marsden(1997)]{KoMa1997c}
W.~S. Koon and J.~E. Marsden.
\newblock The {H}amiltonian and {L}agrangian approaches to the dynamics of
  nonholonomic systems.
\newblock \emph{Reports on Mathematical Physics}, 40\penalty0 (1):\penalty0
  21--62, 1997.

\bibitem[K\"unzle(1969)]{Ku1969}
H.~P. K\"unzle.
\newblock Degenerate {L}agrangean systems.
\newblock \emph{Annales de l'institut Henri Poincar\'e (A)}, 11\penalty0
  (4):\penalty0 393--414, 1969.

\bibitem[Lall and West(2006)]{LaWe2006}
S.~Lall and M.~West.
\newblock Discrete variational {H}amiltonian mechanics.
\newblock \emph{Journal of Physics A: Mathematical and General}, 39\penalty0
  (19):\penalty0 5509--5519, 2006.

\bibitem[Leimkuhler and Reich(2004)]{LeRe2004}
B.~Leimkuhler and S.~Reich.
\newblock \emph{Simulating {H}amiltonian dynamics}, volume~14 of
  \emph{Cambridge Monographs on Applied and Computational Mathematics}.
\newblock Cambridge University Press, Cambridge, 2004.

\bibitem[Leok and Ohsawa(2010)]{LeOh2010}
M.~Leok and T.~Ohsawa.
\newblock Discrete {D}irac structures and implicit discrete {L}agrangian and
  {H}amiltonian systems.
\newblock In \emph{XVIII International Fall Workshop on Geometry and Physics},
  volume 1260, pages 91--102. AIP, 2010.

\bibitem[Leok et~al.()Leok, Ohsawa, and Sosa]{DiracHJ}
M.~Leok, T.~Ohsawa, and D.~Sosa.
\newblock {H}amilton--{J}acobi theory for degenerate {L}agrangian systems with
  constraints.
\newblock \emph{in preparation}.

\bibitem[Lew et~al.(2003)Lew, Marsden, Ortiz, and West]{LeMaOrWe2003}
A.~Lew, J.~E. Marsden, M.~Ortiz, and M.~West.
\newblock An overview of variational integrators.
\newblock In \emph{Finite Element Methods: 1970's and Beyond}. CIMNE, 2003.

\bibitem[Marrero et~al.(2006)Marrero, Mart\'in~de Diego, and
  Mart\'inez]{MaMaMa2006}
J.~C. Marrero, David Mart\'in~de Diego, and E.~Mart\'inez.
\newblock Discrete {L}agrangian and {H}amiltonian mechanics on {L}ie groupoids.
\newblock \emph{Nonlinearity}, 19\penalty0 (6):\penalty0 1313--1348, 2006.

\bibitem[Marsden and Ratiu(1999)]{MaRa1999}
J.~E. Marsden and T.~S. Ratiu.
\newblock \emph{Introduction to Mechanics and Symmetry}.
\newblock Springer, 1999.

\bibitem[Marsden and West(2001)]{MaWe2001}
J.~E. Marsden and M.~West.
\newblock Discrete mechanics and variational integrators.
\newblock \emph{Acta Numerica}, pages 357--514, 2001.

\bibitem[Marsden et~al.(1999)Marsden, Pekarsky, and Shkoller]{MaPeSh1999}
J.~E. Marsden, S.~Pekarsky, and S.~Shkoller.
\newblock Discrete {E}uler--{P}oincar\'e and {L}ie--{P}oisson equations.
\newblock \emph{Nonlinearity}, 12\penalty0 (6):\penalty0 1647--1662, 1999.

\bibitem[McLachlan and Perlmutter(2006)]{McPe2006}
R.~McLachlan and M.~Perlmutter.
\newblock Integrators for nonholonomic mechanical systems.
\newblock \emph{Journal of Nonlinear Science}, 16\penalty0 (4):\penalty0
  283--328, 2006.

\bibitem[Ohsawa and Bloch(2009)]{OhBl2009}
T.~Ohsawa and A.~M. Bloch.
\newblock Nonholonomic {H}amilton--{J}acobi equation and integrability.
\newblock \emph{Journal of Geometric Mechanics}, 1\penalty0 (4):\penalty0
  461--481, 2009.

\bibitem[Ohsawa et~al.(2011)Ohsawa, Fernandez, Bloch, and Zenkov]{OhFeBlZe2011}
T.~Ohsawa, O.~E. Fernandez, A.~M. Bloch, and D.~V. Zenkov.
\newblock Nonholonomic {H}amilton--{J}acobi theory via {C}haplygin
  {H}amiltonization.
\newblock \emph{Journal of Geometry and Physics}, 61\penalty0 (8):\penalty0
  1263--1291, 2011.

\bibitem[Stern(2010)]{St2010}
A.~Stern.
\newblock Discrete {H}amilton--{P}ontryagin mechanics and generating functions
  on {L}ie groupoids.
\newblock \emph{J. Symplectic Geom.}, 8\penalty0 (2):\penalty0 225--238, 2010.

\bibitem[Tulczyjew(1976{\natexlab{a}})]{Tu1976a}
W.~M. Tulczyjew.
\newblock Les sous-vari\'et\'es lagrangiennes et la dynamique hamiltonienne.
\newblock \emph{C. R. Acad. Sc. Paris}, 283:\penalty0 15--18,
  1976{\natexlab{a}}.

\bibitem[Tulczyjew(1976{\natexlab{b}})]{Tu1976b}
W.~M. Tulczyjew.
\newblock Les sous-vari\'et\'es lagrangiennes et la dynamique lagrangienne.
\newblock \emph{C. R. Acad. Sc. Paris}, 283:\penalty0 675--678,
  1976{\natexlab{b}}.

\bibitem[van~der Schaft(1998)]{Sc1998}
A.~J. van~der Schaft.
\newblock Implicit {H}amiltonian systems with symmetry.
\newblock \emph{Reports on Mathematical Physics}, 41\penalty0 (2):\penalty0
  203--221, 1998.

\bibitem[van~der Schaft(2006)]{Sc2006}
A.~J. van~der Schaft.
\newblock Port-{H}amiltonian systems: an introductory survey.
\newblock In \emph{Proceedings of the International Congress of
  Mathematicians}, volume~3, pages 1339--1365, 2006.

\bibitem[van~der Schaft and Maschke(1994)]{ScMa1994}
A.~J. van~der Schaft and B.~M. Maschke.
\newblock On the {H}amiltonian formulation of nonholonomic mechanical systems.
\newblock \emph{Reports on Mathematical Physics}, 34\penalty0 (2):\penalty0
  225--233, 1994.

\bibitem[Vankerschaver et~al.(2010)Vankerschaver, Yoshimura, and
  Marsden]{VaYoMa2010}
J.~Vankerschaver, H.~Yoshimura, and J.~E. Marsden.
\newblock Multi-{D}irac structures and {H}amilton--{P}ontryagin principles for
  {L}agrange--{D}irac field theories.
\newblock \emph{Preprint}, {\tt arXiv:1008.0252}, 2010.

\bibitem[Varadarajan(1984)]{Va1984}
V.~S. Varadarajan.
\newblock \emph{Lie groups, {L}ie algebras, and their representations}.
\newblock Springer, New York, 1984.

\bibitem[Weinstein(1996)]{We1996a}
A.~Weinstein.
\newblock Lagrangian mechanics and groupoids.
\newblock \emph{Fields Inst. Commun.}, 7:\penalty0 207--231, 1996.

\bibitem[Yoshimura and Marsden(2006{\natexlab{a}})]{YoMa2006a}
H.~Yoshimura and J.~E. Marsden.
\newblock {D}irac structures in {L}agrangian mechanics {P}art {I}: Implicit
  {L}agrangian systems.
\newblock \emph{Journal of Geometry and Physics}, 57\penalty0 (1):\penalty0
  133--156, 2006{\natexlab{a}}.

\bibitem[Yoshimura and Marsden(2006{\natexlab{b}})]{YoMa2006b}
H.~Yoshimura and J.~E. Marsden.
\newblock {D}irac structures in {L}agrangian mechanics {P}art {II}: Variational
  structures.
\newblock \emph{Journal of Geometry and Physics}, 57\penalty0 (1):\penalty0
  209--250, 2006{\natexlab{b}}.

\bibitem[Yoshimura and Marsden(2007{\natexlab{a}})]{YoMa2007}
H.~Yoshimura and J.~E. Marsden.
\newblock Reduction of {D}irac structures and the {H}amilton-{P}ontryagin
  principle.
\newblock \emph{Reports on Mathematical Physics}, 60\penalty0 (3):\penalty0
  381--426, 2007{\natexlab{a}}.

\bibitem[Yoshimura and Marsden(2007{\natexlab{b}})]{YoMa2007b}
H.~Yoshimura and J.~E. Marsden.
\newblock {D}irac structures and the {L}egendre transformation for implicit
  {L}agrangian and {H}amiltonian systems.
\newblock In \emph{{L}agrangian and {H}amiltonian Methods for Nonlinear Control
  2006}, pages 233--247, 2007{\natexlab{b}}.

\bibitem[Yoshimura and Marsden(2009)]{YoMa2009}
H.~Yoshimura and J.~E. Marsden.
\newblock Dirac cotangent bundle reduction.
\newblock \emph{Journal of Geometric Mechanics}, 1\penalty0 (1):\penalty0
  87--158, 2009.

\end{thebibliography}
\bibliographystyle{plainnat}

\end{document}